\magnification=1200
\def\A{{\cal A}}
\def\C{{\bf C}}
\def\CC{{\cal C}}
\def\I{{\cal I}}

\def\M{{\cal M}}
\def\N{{\bf N}}
\def\R{{\bf R}}
\def\S{{\cal S}}
\def\T{{\bf T}}
\def\Z{{\bf Z}}
\def\hal{{\vrule height 10pt width 4pt depth 0pt}}

\centerline{Operator algebras associated with the}
\centerline{Klein-Gordon position representation}
\centerline{in relativistic quantum mechanics}
\medskip

\centerline{Nik Weaver\footnote{*}{Partially supported by NSF grant DMS-0070634
\hfill\break
Math Subject Classification numbers: Primary 47L30, 47L15; Secondary 81R15, 47L90, 46L52, 46J99}}
\bigskip
\bigskip

{\narrower{
\noindent \it We initiate a mathematically rigorous study of Klein-Gordon
position operators in single-particle relativistic quantum mechanics.
Although not self-adjoint, these operators have real spectrum
and enjoy a limited form of spectral decomposition. The
associated C*-algebras are identifiable as crossed products. We also
introduce a variety of non self-adjoint operator algebras associated with
the Klein-Gordon position representation; these algebras are commutative
and continuously (but not homeomorphically) embeddable in corresponding
function algebras. Several open problems are indicated.
\bigskip}}
\bigskip

\noindent {\bf 1. Physical background}
\bigskip

Consider a free, spinless, nonrelativistic quantum mechanical particle
in $\R^n$. Its state is represented by a normalized function $\psi$ in
$L^2(\R^n)$, such that the probability of the particle being found in a
region $S \subset \R^n$ is
$$\int_S |\psi(\vec{x})|^2\, d\vec{x}.$$
Similarly, the probability that a measurement of its momentum will return
a value in $S$ is given by
$$\int_S |\phi(\vec{p})|^2\, d\vec{p}$$
where $\phi = \hat{\psi}$ is the Fourier transform of $\psi$. (We adopt
units in which $\hbar = c = 1$ and normalize the Fourier transform to
be unitary from $L^2(\R^n)$ to $L^2(\R^n)$.) The associated nonrelativistic
position operators are the multiplication operators $M_{x_i}$
($1 \leq i \leq n$) on $L^2(\R^n)$.
\medskip

The relativistic picture is more complicated (see, e.g., [20]). The momentum
representation now involves the weighted $L^2$ space $L^2(\R^n, d\vec{p}/E)$
where $E = E(\vec{p}) = \sqrt{|\vec{p}|^2 + m^2}$ is the relativistic
energy, $m$ being the mass of the particle. (We consider throughout
only the simplest case of positive mass and zero spin.) Thus, the probability
of a momentum measurement yielding a value in $S \subset \R^n$ is
$$\int_S |\phi(\vec{p})|^2\, {{d\vec{p}}\over{\sqrt{|\vec{p}|^2 + m^2}}},$$
where $\phi \in L^2(\R^n, d\vec{p}/E)$ is the normalized state of the
particle in this momentum representation. The use of the measure
$d\vec{p}/E$
is dictated by the requirement that Lorentz transformations must induce
unitary transformations of state space; one can regard $\phi$ as being a
function on the Lorentz-invariant positive mass shell
$$\{(\vec{p}, E) \in \R^{n+1}: E > 0, E^2 - |\vec{p}|^2 = m^2\},$$
on which $d\vec{p}/E$ is the essentially unique Lorentz-invariant measure.
\medskip

The problem of determining a valid relativistic position representation has
long been controversial. Those oriented in the direction of mathematical
rigor generally seem to prefer the representation advanced by Newton and
Wigner [15] and Wightman [22]. In this picture the position representation
of the state of a particle with momentum representation
$\phi \in L^2(\R^n, d\vec{p}/E)$ is
$$\psi_{NW} = (\phi/\sqrt{E})\check{\phantom{i}} \in L^2(\R^n),$$
and the position operators are exactly the same as in the nonrelativistic
case. Here $\check{\phantom{a}}$ denotes the inverse Fourier transform. This
representation has the advantage that the map $\phi \mapsto \psi_{NW}$ is
unitary from $L^2(\R^n, d\vec{p}/E)$ to $L^2(\R^n)$. Moreover, Newton,
Wigner, and Wightman established that their position representation was the
only one which satisfied a small set of seemingly undebatable postulates.
\medskip

Nonetheless, most physicists have evidently continued to prefer the
so-called ``Klein-Gordon'' position representation [19], in which the
state is represented as
$$\psi_{KG} = \check{\phi}.$$
(See [13] for a review of other proposed position representations, and
[3] for more recent references.) Among
the difficulties associated with this representation are the fact that
$\check{\phi}$ need not belong to $L^2(\R^n)$, or even to be well-defined
except as a distribution, and the failure of the corresponding position
operators (these will be described in the next section) to be self-adjoint.
\medskip

A third oddity of the Klein-Gordon position representation, and the cause
of its failure to satisfy the Newton-Wigner-Wightman postulates, is the
fact that states localized in disjoint spatial regions are in general not
orthogonal. At first sight this seems to violate relativistic causality.
However, the mere nonorthogonality of distant spatially separated states
does not actually imply any mechanism for superluminal signalling.
Moreover, a frustrating and seemingly unavoidable feature of essentially
any relativistic notion of position for single particles --- including
both the Newton-Wigner and the Klein-Gordon proposals --- is the
instantaneous spreading throughout space of any free particle which is
localized in a bounded spatial region at some time [9]. In particular,
the state of a free particle localized in some bounded region at some
time in the Newton-Wigner sense will, in general, not be orthogonal to a
state which is localized in an arbitrarily distant region at a slightly
later time. This fact calls into question the necessity of requiring
orthogonality at a fixed time and makes the nonorthogonality exhibited
by Klein-Gordon localized particles seem less objectionable in comparison
to the Newton-Wigner picture.
\medskip

Philips [16] took advantage of this point by proposing an alternative
set of postulates for localization, which are satisfied by the Klein-Gordon
representation but not the Newton-Wigner representation and which replace
the orthogonality requirement mentioned above by a requirement of ``Lorentz
invariance.'' Informally, this condition asserts that a state which appears
as a delta function (i.e., is localized at a single spatial point) at some
time according to one inertial observer must appear as a delta function at
the same space-time event according to all inertial observers. Now since
delta states are non-physical, it seems unlikely that this postulate has
any direct physical content. In particular, one cannot reformulate the
Lorentz invariance condition in terms of states localized in arbitrarily
small but finite spatial regions: that is, one cannot require
that states supported sufficiently near a given space-time
event according to one observer must be supported near that event
according to all observers. This is because the phenomenon of
instantaneous propagation mentioned above implies that a state localized
in a small region at some time according to one inertial observer will in
fact not be localized in any bounded region at any time according to other
observers, regardless of the position representation used. It therefore
appears that, contrary to a claim often made in the literature, the
requirement of Lorentz invariance is not actually necessary on physical
grounds, nor even directly physically interpretable.
\medskip

Thus, the fact that localized states instantaneously propagate throughout
space calls into question, in different ways, the physical justification
for both the Newton-Wigner and the Klein-Gordon position representations.
Analyzing the localization problem in terms of field theory also fails to
resolve the conflict: although one has mathematically rigorous models of
free fields which identify the observables that can be measured in any
given spatial region ([8]; see [21] for an elementary treatment), it
turns out that for any spatial region $S$, the entire one-particle
Hilbert space of a free field is generated by those one-particle states
which can be detected with certainty by field observables localized in
$S$. This phenomenon suggests that definite physical localization of
one-particle states is simply not possible, so that the debate between
rival one-particle localization schemes is operationally meaningless.
\medskip

Of course, this negative conclusion does not deny that various position
representations may be theoretically useful in different ways. However,
if the issue is one of theoretical utility, and this apparently is the
only real issue, we must face the uncomfortable fact that the Newton-Wigner
notion of position, despite its manifest advantages from the mathematical
point of view, has not been taken up by physicists with much enthusiasm.
In fact the vast majority of works on relativistic quantum mechanics
adopt the Klein-Gordon representation without comment, and the bulk of
those which do discuss position or localization explicitly also favor
this representation (see, e.g., [2], [7], [10], [11], [12], [14], [18]).
Even a brief glance at any standard reference on relativistic quantum
mechanics will show how essential the Klein-Gordon representation is,
and how peripheral the Newton-Wigner representation, to its usual treatment.
\medskip

In this light, it is rather surprising that the Klein-Gordon position
representation for one-particle states has, up to now, not been explored
from a serious mathematical point of view. The goal of this paper is to
begin to study the Klein-Gordon representation in a mathematically rigorous
manner, and in particular to place it in an operator algebra context. As
will be seen, many interesting questions emerge, several of which we have
been unable to answer. Moreover, especially in later sections, the reader
will find it easy to pose natural questions not addressed here.
\bigskip
\bigskip

\noindent {\bf 2. Spectral subspaces}
\bigskip

We will mainly work in the Fourier transformed (i.e., momentum)
picture. As explained in Section 1, this involves the Hilbert space
$L^2(\R^n, d\vec{p}/E)$. Now the Klein-Gordon position operators are,
heuristically, coordinate multiplication operators in the untransformed
picture, and they can therefore be rigorously defined as differentiation
operators in the transformed picture. Similarly, translations and
convolutions in the transformed picture can be viewed as different
kinds of multiplication operators in the untransformed picture.
\bigskip

\noindent {\bf Definition 1.}
\smallskip

\noindent {\bf (a)} The {\it relativistic momentum Hilbert space}
is the space $L^2(\R^n, d\vec{p}/E)$, where $E(\vec{p}) =
\sqrt{|\vec{p}|^2 + m^2}$ and $m > 0$ is fixed. For simplicity we
will write $L^2_r(\R^n) = L^2(\R^n, d\vec{p}/E)$. The norm of
$\phi \in L^2_r(\R^n)$ will be denoted
$$\|\phi\|_r = \left(\int |\phi(\vec{p})|^2
{{d\vec{p}}\over{E(\vec{p})}}\right)^{1/2}.$$
\medskip

\noindent {\bf (b)} The {\it Klein-Gordon position operators} are
the unbounded operators $Q_i$ ($1 \leq i \leq n$) on $L^2_r(\R^n)$
defined by
$$Q_i \phi = \iota{{\partial\phi}\over{\partial p_i}},$$
with domain $D(Q_i)$ the set of $\phi \in L^2_r(\R^n)$ whose restrictions
to lines parallel to the $i$th coordinate axis are almost all locally
absolutely continuous, and which satisfy $\partial\phi/\partial p_i
\in L^2_r(\R^n)$.
\medskip

\noindent {\bf (c)} For $\vec{a} \in \R^n$ let $T_{\vec{a}}$ be the
translation operator on $L^2_r(\R^n)$ defined by
$T_{\vec{a}} \phi(\vec{p}) = \phi(\vec{p} - \vec{a})$.
\medskip

\noindent {\bf (d)} Let $\S(\R^n)$ denote the class of Schwartz functions
on $\R^n$, and for $f \in \S(\R^n)$ let $K_f$ be the convolution operator
(with respect to Lebesgue measure) acting on $L^2_r(\R^n)$ defined by
$$K_f \phi(\vec{p}) = f * \phi(\vec{p})
= \int_{\R^n} f(\vec{p}\,{}')\phi(\vec{p} - \vec{p}\,{}')\, d\vec{p}\,{}'.$$

In line with the previous comment, $T_{\vec{a}}$ is to be thought of
heuristically as multiplication by $e^{\iota\vec{a}\cdot\vec{x}}$, and
$K_f$ as multiplication by $\check{f}$, in the position representation
(untransformed picture).
\bigskip

\noindent {\bf Remarks 2.}
\smallskip

\noindent {\bf (a)} The usual Hilbert space $L^2(\R^n)$ is densely
contained in $L^2_r(\R^n)$. Moreover, the set of compactly supported
functions in $L^2(\R^n)$ equals the set of compactly supported functions
in $L^2_r(\R^n)$.
\medskip

\noindent {\bf (b)} An equivalent definition of $D(Q_i)$ is:
$\phi \in D(Q_i)$ if for every $N > 0$ the restriction $\phi|_{[-N,N]^n}$
belongs to the maximal domain of $\iota\partial/\partial p_i$ on
$L^2([-N,N]^n)$, and $\partial\phi/\partial p_i \in L^2_r(\R^n)$.
\medskip

\noindent {\bf (c)} The adjoint of $Q_i$ is the operator $M_EQ_iM_{1/E}$,
i.e., the operator $\phi \mapsto E\cdot(\iota\partial/\partial p_i)(\phi/E)$.
Thus, $Q_i$ is not self-adjoint, nor even normal.
\medskip

\noindent {\bf (d)} A simple calculation shows that the operators
$T_{\vec{a}}$ ($\vec{a} \in \R^n$) are bounded, and the $K_f$
($f \in \S(\R^n)$) are bounded by [5, Theorem II.1.6] (cf.\ Lemma 36).
Note that the $T_{\vec{a}}$ are not unitary.
\bigskip

Although the $Q_i$ are not normal, it is interesting to note that their
spectra are real, just as in the nonrelativistic case.
\bigskip

\noindent {\bf Proposition 3.} {\it For each $i$, the operator $Q_i$
is closed and its spectrum is $\R$.}
\medskip

\noindent {\it Proof.} Both claims are fairly standard. To verify closure,
suppose $\phi_k \to \phi$ and $Q_i\phi_k \to \psi$ in $L^2_r(\R^n)$. If
$S = [-N,N]^n \subset \R^n$ then $L^2(S, d\vec{p}/E) \cong L^2(S)$, so
$\phi_k|_S \to \phi|_S$ and $(Q_i\phi_k)|_S \to \psi|_S$ in $L^2(S)$.
Letting $A$ be the maximal version of $\iota\partial \phi/\partial p_i$ on
$L^2(S)$, we have $\phi_k|_S \in D(A)$ and $A(\phi_k|_S) = (Q_i\phi_k)|_S$
for all $k$. Closure of $A$ implies that $\phi|_S \in D(A)$ and
$\iota\partial \phi/\partial p_i = \psi$ on $S$. It follows that
$\phi \in D(Q_i)$ and $Q_i\phi = \psi$, so $Q_i$ is closed.
\medskip

Let $z \in \C$ and suppose ${\rm Im}\, z > 0$. Then the operator
$R_z$ defined by
$$R_z\phi(\vec{p}) =
\iota\int_0^\infty e^{\iota tz}\phi(\vec{p} + t\vec{e}_i)\, dt,$$
where $\vec{e}_i$ is the canonical basis vector, is a bounded
inverse of $Q_i - zI$. A similar expression accomplishes the same result
when ${\rm Im}\, z < 0$. Thus the spectrum of $Q_i$ is contained in $\R$.
Conversely, for any $a \in \R$ the functions
$$f_k(\vec{p}) = e^{\iota ap_i}e^{-|\vec{p}|^2/k}$$
are approximate eigenvectors for $Q_i$ with approximate eigenvalue $a$.
This implies that $(Q_i - aI)^{-1}$ is not bounded. So the spectrum of
$Q_i$ equals $\R$.\hfill\hal
\bigskip

In fact the $Q_i$ even have simultaneous approximate eigenvectors, for
instance the functions $e^{\iota\vec{a}\cdot\vec{p}}e^{-|\vec{p}|^2/k}$.
This raises the question, to what extent do the $Q_i$ support a spectral
decomposition of $L^2_r(\R^n)$? The following seems a reasonable notion
of spectral subspaces associated to the Klein-Gordon position operators.
\bigskip

\noindent {\bf Definition 4.} For any positive measure set $S \subset
\R^n$, let $H_S \subset L^2_r(\R^n)$ be the closure (with respect to
$\|\cdot\|_r$) of the set $\{\hat{\psi}: \psi \in L^2(S)\}$.
\bigskip

\noindent {\bf Proposition 5.} {\it Let $S \subset \R^n$ be a bounded
measurable set. Then $H_S$ is contained in the joint domain of the
$Q_i$, and the distributional equation $Q_i\phi =
(x_i\check{\phi})\hat{\phantom{i}}$ holds for all $\phi \in H_S$.}
\medskip

\noindent {\it Proof.} Fix $1 \leq i \leq n$ and $f \in \S(\R^n)$ such
that $f(\vec{x}) = x_i$ on $S$. Then let $\phi \in H_S$ and find a
sequence $(\psi_k)$ in $L^2(S)$ such that $\hat{\psi}_k \to \phi$ in
$L^2_r(\R^n)$. We have $f\psi_k = x_i\psi_k$ for all $k$, and hence
$K_{\hat{f}} \hat{\psi}_k = Q_i\hat{\psi}_k$ for all $k$. This implies
that the sequence $(Q_i\hat{\psi}_k)$ converges in $L^2_r(\R^n)$ (since
$K_{\hat{f}}$ is bounded) and therefore, since $Q_i$ is closed, that
$\phi \in D(Q_i)$. The equation $Q_i\phi =
(x_i\check{\phi})\hat{\phantom{i}}$ holds by continuity.\hfill\hal
\bigskip

Now disjointness of $S$ and $T$ does not imply that $H_S$
and $H_T$ are orthogonal; this corresponds to the nonorthogonality
of disjointly supported states mentioned in Section 1. However, we
can still hope that $H_S \cap H_T = \{0\}$ when $S$ and $T$ are
disjoint. This is at least true to the following extent:
\bigskip

\noindent {\bf Proposition 6.} {\it Let $S, T \subset \R^n$ be
disjoint closed sets. Then $H_S \cap H_T = \{0\}$.}
\medskip

\noindent {\it Proof.} Let $\phi \in H_S \cap H_T$ and suppose
$\phi \neq 0$. Fix $\vec{x}$ in the support of the distribution
$\check{\phi}$. Since $S$ and $T$ are disjoint,
without loss of generality we can assume $\vec{x} \not\in S$.
Let $f$ be a smooth, compactly supported function which is
constantly zero on $S$ and satisfies $f(\vec{x}) \neq 0$. Then
$f\check{\phi} \neq 0$, so $\hat{f} * \phi \neq 0$ by the definition of
the convolution of a function with a distribution. But since
$\phi \in H_S$ we can find a sequence $(\psi_k)$ in $L^2(S)$
such that $\hat{\psi}_k \to \phi$ in $L^2_r(\R^n)$. Then
$f\psi_k = 0$ for all $k$, which implies
$$K_{\hat{f}}\hat{\psi}_k = (f\psi_k)\hat{\phantom{i}} = 0$$
for all $k$, which implies $K_{\hat{f}}\phi = 0$ by continuity. This
is a contradiction, and we conclude that $\phi = 0$.\hfill\hal
\bigskip

The technique of the preceding proof can be adapted to the following
situation as well.
\bigskip

\noindent {\bf Proposition 7.} {\it Let $S, T \subset \R^n$ be disjoint
measurable sets and suppose the boundary of $S$ is piecewise linear. Then
$H_S \cap H_T = \{0\}$.}
\medskip

\noindent {\it Proof.} Let $\phi \in H_S \cap H_T$ and suppose
$\phi \neq 0$. The argument in the proof of Proposition 6 shows
that the support of $\check{\phi}$ must be contained in both
the closure of $S$ and the closure of $T$, and hence it must be
contained in the boundary of $S$. Since this boundary is piecewise
linear, we can find a smooth
function $f$ with compact support such that $f\check{\phi}$ is
nonzero and is supported in a hyperplane. But the Fourier transform
of $f\check{\phi}$ then has constant modulus in directions perpendicular
to the hyperplane, so $\hat{f}*\phi = K_{\hat{f}}\phi$ cannot belong
to $L^2_r(\R^n)$, contradicting boundedness of $K_{\hat{f}}$. We
conclude that $\phi = 0$.\hfill\hal
\bigskip

\noindent {\bf Problem 8.} Does $H_S \cap H_T = \{0\}$ hold for any
disjoint positive measure sets $S, T \subset \R^n$? Does it hold if
$S$ has smooth boundary?
\bigskip

We note that the answer to Problem 8 depends on the rate of decrease
of the weight function $1/E(\vec{p})$. If we took $E(\vec{p}) \equiv 1$,
then $L^2_r(\R^n)$
would equal $L^2(\R^n)$ and we would certainly have $H_S \cap H_T = \{0\}$
whenever $S$ and $T$ were disjoint. On the other hand, if we took
$E(p) = 1 + p^2$ in the case $n = 1$, then the function $\phi(p) \equiv 1$
would belong to $L^2_r(\R)$ and would be in $H_S \cap H_T$ for
$S = (0,1)$ and $T = (-1, 0)$. This shows that even Proposition 7 does
not hold for weights of more rapid decrease.
\bigskip
\bigskip

\noindent {\bf 3. C*- and von Neumann algebras}
\bigskip

In this section we identify the structure of some self-adjoint operator
algebras affiliated with the Klein-Gordon position operators. Although
the non self-adjoint operator algebras to be discussed in later sections
seem more interesting, we deal with the self-adjoint case first since
these algebras are easier.
\medskip

Specifically, the algebras we consider in this section are the
C*-algebras
$$C^*(T_{\vec{a}}: \vec{a} \in \R^n)$$
and
$$C^*(K_f: f \in \S(\R^n))$$
generated by the operators $T_{\vec{a}}$ and $K_f$, and the von Neumann
algebra
$$W^*(T_{\vec{a}}: \vec{a} \in \R^n) = W^*(K_f: f \in \S(\R^n))$$
generated by either set. (It is straightforward to verify that every
$T_{\vec{a}}$ can be approximated in the strong operator topology by
operators of the form $K_f$, and that every $K_f$ can be strong operator
approximated by finite linear combinations of operators of the form
$T_{\vec{a}}$ --- just approximate the integral $K_f =
\int f(\vec{a})T_{\vec{a}}\, d\vec{a}$ by Riemann sums ---
so that the two definitions of the von Neumann algebra do coincide.)
\medskip

In order to analyze these algebras, it is helpful to pass to a different
Hilbert space. The desired transformation is this:
\bigskip

\noindent {\bf Definition 9.} Regard the multiplication operator
$M_{\sqrt{E}}$ as a unitary map from the usual Hilbert space $L^2(\R^n)$
to the relativistic space $L^2_r(\R^n)$. Then define operators
$\tilde{T}_{\vec{a}}$ and $\tilde{K}_f$ on $L^2(\R^n)$ by
$$\tilde{T}_{\vec{a}} = M_{\sqrt{E}}^{-1}T_{\vec{a}}M_{\sqrt{E}}
\qquad{\rm and}\qquad \tilde{K}_f = M_{\sqrt{E}}^{-1}K_fM_{\sqrt{E}}$$
for $\vec{a} \in \R^n$ and $f \in \S(\R^n)$.
\bigskip

We first show that the von Neumann algebra generated by the $T_{\vec{a}}$
(or equivalently, the $K_f$) is trivial.
\bigskip

\noindent {\bf Theorem 10.} {\it
$W^*(\tilde{T}_{\vec{a}}: \vec{a} \in \R^n) = B(L^2(\R^n))$.}
\medskip

\noindent {\it Proof.} Let $\M = W^*(\tilde{T}_{\vec{a}}: \vec{a} \in \R^n)
\subset B(L^2(\R^n))$.
Simple calculations show that
$$\tilde{T}_{\vec{a}}\phi(\vec{p})
= \sqrt{{E(\vec{p} - \vec{a})}\over{E(\vec{p})}} \phi(\vec{p} - \vec{a})$$
and
$$\tilde{T}^*_{\vec{a}}\phi(\vec{p})
= \sqrt{{E(\vec{p})}\over{E(\vec{p} + \vec{a})}} \phi(\vec{p} + \vec{a})$$
for all $\phi \in L^2(\R^n)$. Therefore the operator
$\tilde{T}^*_{\vec{a}}\tilde{T}_{\vec{a}}$ is multiplication by
the function ${{E(\vec{p})}/{E(\vec{p} + \vec{a})}}$. These functions
separate points in $\R^n$, and therefore $\M$ must contain all
bounded multiplication operators.
In particular, it contains multiplication by the function
$\sqrt{E(\vec{p})/E(\vec{p} - \vec{a})}$ for every $\vec{a}$; taking
the product of this operator with $\tilde{T}_{\vec{a}}$, we deduce
that every translation operator on $L^2(\R^n)$ belongs to $\M$.
\medskip

Now since $\M$ contains every multiplication operator,
any operator in the commutant of $\M$ must itself be a
multiplication operator. But $\M$ also contains all translations,
and the only multiplication operators which commute with all
translations are the scalar multiples of the identity.
Thus $\M' = \C\cdot I$ and $\M = \M'' = B(L^2(\R^n))$.\hfill\hal
\bigskip

Next we characterize the structure of the two C*-algebras mentioned above.
\bigskip

\noindent {\bf Theorem 11.} {\it $C^*(\tilde{T}_{\vec{a}}: \vec{a} \in \R^n)$
is the C*-algebra generated by the translations on $L^2(\R^n)$ together with
the multiplication operators by continuous functions vanishing at infinity. It
is $*$-isomorphic to the crossed product $C_0(\R^n)^+ \times_\lambda \R^n_d$,
where $C_0(\R^n)^+$ is the unitization of $C_0(\R^n)$, $\R^n_d$ is $\R^n$
with the discrete topology, and $\lambda$ is the action of $\R^n_d$ on
$C_0(\R^n)^+$ by translations.}
\medskip

\noindent {\it Proof.} Let $\A_1 = C^*(\tilde{T}_{\vec{a}}: \vec{a} \in \R^n)$.
As in the proof of Theorem 10, $\A_1$ contains multiplication by the
function ${{E(\vec{p})}/{E(\vec{p} + \vec{a})}}$, for every $\vec{a} \in \R^n$.
These functions separate points and they all converge to 1 at infinity
(and are not all 1 at any other point), so
they generate the C*-algebra of multiplications by continuous functions
on $\R^n$ which are continuously extendible to the one-point compactification
of $\R^n$. In particular, $\A_1$ contains multiplication by any function in
$C_0(\R^n)$; also, as in the proof of Theorem 10 we can deduce that $\A_1$
contains all translation operators.
Thus one containment of the first assertion of the theorem
holds. The reverse containment is easy: any C*-algebra which contains
the translations must contain the identity operator, and if it also
contains multiplications by functions in $C_0(\R^n)$ then it contains
multiplications by all continuous functions which converge at infinity, in
particular by the functions $\sqrt{E(\vec{p} - \vec{a})/E(\vec{p})}$.
From the expression for $\tilde{T}_{\vec{a}}$ given in the proof of
Theorem 10, we conclude that the C*-algebra generated by $C_0(\R^n)$
and the translations contains $\A_1$. This completes the proof of the
first assertion of the theorem.
\medskip

For the second assertion, recall that the crossed product
$C_0(\R^n)^+ \times_\gamma \R^n_d$ is the universal C*-algebra generated
by a copy of $C_0(\R^n)^+$ and a (not continuous) one-parameter group of
unitaries $U_{\vec{a}}$ ($\vec{a} \in \R^n$) which satisfy the commutation
relations $U_{\vec{a}}f = f_{\vec{a}} U_{\vec{a}}$, where $f_{\vec{a}}$
is the translation of $f \in C_0(\R^n)^+$ by $\vec{a}$. Now $\A_1$ is,
by the preceding
paragraph, generated by a copy of $C_0(\R^n)^+$ and the (unitary)
translations on $L^2(\R^n)$, which satisfy the same commutation relations.
Thus $\A_1$ is naturally a $*$-homomorphic image of the crossed product
C*-algebra. To show isomorphism, we must verify that the kernel of this
$*$-homomorphism is zero.
\medskip

Let $\I$ be this kernel. Consider the action $\alpha$ of $\R^n$ on
$C_0(\R^n)^+ \times_\gamma \R^n_d$ defined by
$\alpha_{\vec{a}}(fU_{\vec{b}}) = e^{\iota\vec{a}\cdot\vec{b}}fU_{\vec{b}}$.
Every element of the crossed product is almost periodic for this action,
so taking mean-value integrals yields
a mean value map $m$ from the crossed product to the fixed point
algebra for $\alpha$. But the mean value map is the identity on
$C_0(\R^n)^+$ and annihilates all elements of the form $fU_{\vec{b}}$
for $\vec{b} \neq 0$. Thus $C_0(\R^n)^+$ is the fixed point algebra.
Now suppose $\I \neq \{0\}$. Then $\I$ contains a positive operator $A$,
and so $m(A) > 0$ by a standard property of the mean value map. Let
$f \in C_0(\R^n)$ be compactly supported and satisfy $f \geq 0$ and
$f^2m(A) \neq 0$.
\medskip

We claim that $\I$ contains $m(fAf) = f^2\cdot m(A)$. It will suffice to
verify that $\I$ contains $\alpha_{\vec{a}}(fAf)$ for all $\vec{a}$. Find
$g_{\vec{a}} \in C_0(\R^n)$ such that $g_{\vec{a}}f =
e^{\iota\vec{a}\cdot\vec{p}}f$. Then $\alpha_{\vec{a}}(fAf) =
g_{\vec{a}}fAf\bar{g}_{\vec{a}} \in \I$ since
$\I$ is an ideal. Thus $\I$ contains $f^2m(A)$, which is a nonzero element
of $C_0(\R^n)^+$. But the quotient map from $C_0(\R^n)^+ \times_\gamma \R^n_d$
to $\A_1$ is clearly isometric on $C_0(\R^n)^+$, a contradiction. We conclude
that $\I$ is zero and the quotient map is a $*$-isomorphism.\hfill\hal
\bigskip

For any operator $A \in B(L^2(\R^n))$, let $\widehat{A}$ denote its
conjugation by the Fourier transform; that is, $\widehat{A}\phi =
(A\check{\phi})\hat{\phantom{i}}$. Also write
$\widehat{\A} = \{\widehat{A}: A \in \A\}$ for any subalgebra $\A$ of
$B(L^2(\R^n))$.
\bigskip

\noindent {\bf Theorem 12.} {\it $C^*(\tilde{K}_f: f \in \S(\R^n))
= \widehat{C_0(\R^n)} + K(L^2(\R^n))$, regarding $C_0(\R^n)$ as acting on
$L^2(\R^n)$ by multiplication. It is $*$-isomorphic to the crossed product
$C_0(\R^n)^+ \times_\lambda \R^n$, where $C_0(\R^n)^+$ is the unitization
of $C_0(\R^n)$ and $\lambda$ is the action of $\R^n$ on $C_0(\R^n)^+$ by
translations (and the group $\R^n$ has its usual topology).}
\medskip

\noindent {\it Proof.} Let $\A_2 = C^*(\tilde{K}_f: f \in \S(\R^n))$.
A straightforward computation shows that
$$\tilde{K}_f \phi(\vec{p}) =
\int f(\vec{p}\,{}') \sqrt{{E(\vec{p} - \vec{p}\,{}')}\over{E(\vec{p})}}
\phi(\vec{p} - \vec{p}\,{}')\, d\vec{p}\,{}'$$
for all $\phi \in L^2(\R^n)$, and its adjoint satisfies
$$\tilde{K}^*_f \phi(\vec{p}) =
\int f^*(\vec{p}\,{}')\sqrt{{E(\vec{p})}\over{E(\vec{p} - \vec{p}\,{}')}}
\phi(\vec{p} - \vec{p}\,{}')\, d\vec{p}\,{}'$$
where $f^*(\vec{p}) = \overline{f(-\vec{p})}$. Thus, for any $f \in
\S(\R^n)$ we have
$$(\tilde{K}_f - \tilde{K}^*_{f^*}) \phi(\vec{p}) =
\int f(\vec{p}\,{}')\left[\sqrt{{E(\vec{p} - \vec{p}\,{}')}\over{E(\vec{p})}} -
\sqrt{{E(\vec{p})}\over{E(\vec{p} - \vec{p}\,{}')}}\,\right]
\phi(\vec{p} - \vec{p}\,{}')\, d\vec{p}\,{}'.$$
This is an integral operator on $L^2(\R^n)$ with kernel in $C_0(\R^{2n})$.
If the kernel function is multiplied by the characteristic function of the
ball of radius $R$ about the origin in $\R^{2n}$, then the corresponding
integral operator will be compact, and these truncated operators converge
to the original operator in norm by an estimate using [5, Theorem II.1.6].
Hence $\tilde{K}_f - \tilde{K}^*_{f^*}$ is compact (and nonzero if $f$ is
nonzero). By Theorem 10, $\A_2$ acts irreducibly on $L^2(\R^n)$. As we have
just shown that $\A_2$ nontrivially intersects the compact operators, it
follows that $\A_2$ contains $K(L^2(\R^n))$ [1, Corollary 2, p.\ 18]. 
\medskip

Now for any $f \in \S(\R^n)$ the operator $A_f \in B(L^2(\R^n))$ defined by
$$A_f\phi(\vec{p}) =
\int f(\vec{p}\,{}') \left[ \sqrt{{E(\vec{p} - \vec{p}\,{}')}\over{E(\vec{p})}}
- 1\right]\phi(\vec{p} - \vec{p}\,{}')\, d\vec{p}\,{}'$$
is another integral operator which is compact by a similar argument to the
one indicated above, and hence it belongs to $\A_2$. It follows that the
difference $K_f - A_f$, which is simply convolution by $f$, also belongs
to $\A_2$ and that $\A_2$ is generated by $K(L^2(\R^n))$ together with the
convolution operators by $f \in \S(\R^n)$. The C*-algebra generated by the
latter operators is just $\widehat{C_0(\R^n)}$, so we conclude that
$\A_2 = \widehat{C_0(\R^n)} + K(L^2(\R^n))$. (Note that we do not have
to take the closure; it is a basic exercise in C*-algebra theory to show
that this sum is already closed.)
\medskip

The crossed product C*-algebra $C_0(\R^n) \times_\lambda \R^n$ is
$*$-isomorphic to $K(L^2(\R^n))$ [6, Theorem 8.4.3]. Thus the
larger crossed product $C_0(\R^n)^+ \times_\lambda \R^n$ is generated
by $K(L^2(\R^n))$ together with $C^*(\R^n) \cong C_0(\R^n)$. This
shows that the latter crossed product is also isomorphic to
$\widehat{C_0(\R^n)} + K(L^2(\R^n))$. (The isomorphism is natural;
it is implemented by the following map. For any $f \in \S(\R^n)$ with
compact support, define a function $\omega_f \in C_c(\R^n, C_0(\R^n)^+)$ by
$$\omega_f(\vec{a}) =
f(\vec{a})\sqrt{{E(\vec{p} - \vec{a})}\over{E(\vec{p})}}.$$
Regarding $C_0(\R^n)^+ \times_\lambda \R^n$ as a completion of
$C_c(\R^n, C_0(\R^n)^+)$, the map $\tilde{K}_f \mapsto \omega_f$
extends to a $*$-isomorphism from $\A_2$ onto
$C_0(\R^n)^+ \times_\lambda \R^n$.)\hfill\hal
\bigskip
\bigskip

\noindent {\bf 4. A non self-adjoint analog of $L^\infty(\R^n)$}
\bigskip

We mentioned at the start of the last section that the self-adjoint
algebras affiliated with the Klein-Gordon position operators seem
less interesting than the non self-adjoint algebras. One
reason for this is that the non self-adjoint algebras play a
functional calculus role that the self-adjoint algebras do not.
Intuitively, the operators $T_{\vec{a}}$ and $K_f$ can be regarded
as arising from the $Q_i$ by the equations $T_{\vec{a}} =
e^{\iota\vec{a}\cdot\vec{Q}}$ and $K_f = \check{f}(\vec{Q})$, where
$\vec{Q} = (Q_1, \ldots, Q_n)$. Their adjoints have no such
interpretation. Thus, various non self-adjoint algebras can be
seen as the set of all operators arising from the $Q_i$ via functional
calculi involving various classes of functions. We begin in this section
with the largest (bounded) class of interest, and discuss other classes
in the next section.
\bigskip

\noindent {\bf Definition 13.}
\smallskip

\noindent {\bf (a)} Denote by $L^\infty_r(\R^n)$ the commutant of the
set $\{T_{\vec{a}}: \vec{a} \in \R^n\}$, or equivalently, of the set
$\{K_f: f \in \S(\R^n)\}$.
\medskip

\noindent {\bf (b)} Let $\gamma$ be the action of $\R^n$ on
$L^\infty_r(\R^n)$ defined by
$$\gamma_{\vec{a}}(A) =
M_{e^{-\iota\vec{a}\cdot\vec{p}}}AM_{e^{\iota\vec{a}\cdot\vec{p}}}.$$

As we noted at the beginning of Section 3, the $K_f$ and linear
combinations of the $T_{\vec{a}}$ mutually strong operator approximate
each other. Thus the commutants of the two sets are the same.
\medskip

The set ${\rm span}\{T_{\vec{a}}: \vec{a} \in \R^n\}$ is
stable under $\gamma$, and it follows that $L^\infty_r(\R^n)$ is also
stable under $\gamma$. The action $\gamma$ is, in the Fourier transform
picture, analogous to the action of $\R^n$ on $L^\infty(\R^n)$ by
translations. It will play a basic role in the sequel.
\bigskip

\noindent {\bf Remarks 14.}
\smallskip
\noindent {\bf (a)} The action $\gamma$ on $L^\infty_r(\R^n)$ arises via
congujation by a unitary action of the group of translations of $\R^n$
on $L^2_r(\R^n)$. The latter extends to a natural representation of the 
$(n+1)$-dimensional Poincar\'e group on $L^2_r(\R^n)$. To see
this, regard elements of $L^2_r(\R^n)$ as functions on the positive mass
shell $\{(\vec{p}, E) \in \R^{n+1}: E > 0, E^2 - |\vec{p}|^2 = m^2\} \subset
\R^{n+1}$. Then Lorentz transformations act on $L^2_r(\R^n)$ by composition,
and these actions are unitary by Lorentz invariance of the measure
$d\vec{p}/E$. The translation of $\R^{n+1}$ by $(\vec{a}, a_{n+1})
\in \R^{n+1}$ acts on $L^2_r(\R^n)$ via the multiplication operator
$M_{e^{\iota(\vec{a}\cdot\vec{p} + a_{n+1}E)}}$.
\medskip

\noindent {\bf (b)} This action of the Poincar\'e group on $L^2_r(\R^n)$
does not induce an action on $L^\infty_r(\R^n)$, since the latter is, in
general, not stable under conjugation by the above unitaries. However, the
group of isometries of $\R^n$, regarded as a subgroup of the Poincar\'e
group, does act on $L^\infty_r(\R^n)$ by isometries in this way.
\medskip

\noindent {\bf (c)} The group $GL(n,\R)$ also acts boundedly on
$L^2_r(\R^n)$ by composition, and conjugation by these operators
yields an action of $GL(n,\R)$ on $L^\infty_r(\R^n)$ by
(non-isometric) isomorphisms.
\bigskip

We begin with a simple fact about $L^\infty_r(\R^n)$.
\bigskip

\noindent {\bf Proposition 15.} {\it $L^\infty_r(\R^n)$ is a
dual space.}
\medskip

\noindent {\it Proof.} It is easy to check that $L^\infty_r(\R^n)$ is a
weak operator closed linear subspace of $B(L^2_r(\R^n))$. This
automatically implies that $L^\infty_r(\R^n)$ is weak* (ultraweakly) closed.
A weak* closed subspace of a dual space is always a dual space.\hfill\hal
\bigskip

We proceed immediately to our main result on $L^\infty_r(\R^n)$,
which states that it consists of convolutions by distributions
which are Fourier transforms of functions in $L^\infty(\R^n)$.
Thus, $L^\infty_r(\R^n)$ can be (non-homeomorphically) embedded
as a subalgebra of $L^\infty(\R^n)$.
\medskip

Recall the notation $\widehat{\A}$ defined just before Theorem 12.
\bigskip

\noindent {\bf Theorem 16.} {\it Regarding $L^2(\R^n)$ as a subset of
$L^2_r(\R^n)$, for every $A \in L^\infty_r(\R^n)$ we have
$A(L^2(\R^n)) \subset L^2(\R^n)$. The restriction $A|_{L^2(\R^n)}$ is
a bounded operator on $L^2(\R^n)$, and the restriction map
$\hat{\sigma}: A \mapsto A|_{L^2(\R^n)}$ is an injective, contractive,
weak*-continuous algebra homomorphism from $L^\infty_r(\R^n)$ into
$\widehat{L^\infty(\R^n)}$.}
\medskip

\noindent {\it Proof.} Let $A \in L^\infty_r(\R^n)$ and let
$\phi \in L^2_r(\R^n)$ be nonzero and
have compact support (so that $\phi \in L^2(\R^n)$
also). Let $T_t$ be the operator of translation by $t$ units in the first
coordinate. Then writing $\|\cdot\|$ for the norm in $L^2(\R^n)$, we have
$$\lim_{t \to \infty} \sqrt{t}\|T_t\phi\|_r = \|\phi\|$$
(using the fact that $\phi$ has compact support), and
$$\liminf_{t \to \infty} \sqrt{t}\|T_tA\phi\|_r \geq \|A\phi\|$$
(even if, hypothetically, $\|A\phi\|$ is infinite), so that
$$\liminf_{t \to \infty} {{\|AT_t\phi\|_r}\over{\|T_t\phi\|_r}}
= \liminf_{t \to \infty} {{\|T_tA\phi\|_r}\over{\|T_t\phi\|_r}}
\geq {{\|A\phi\|}\over{\|\phi\|}}.$$
This shows that $\|A\phi\|/\|\phi\| \leq \|A\|$ (the norm of $A$ as an
element of $L^\infty_r(\R^n)$) for all compactly supported, nonzero $\phi$,
but we cannot immediately deduce that $A|_{L^2(\R^n)}$ is bounded on
$L^2(\R^n)$. For arbitrary nonzero $\phi$ we must make the following
argument. Let
$\phi \in L^2(\R^n)$ and find a sequence $(\phi_k)$ of compactly supported
functions which converge to $\phi$ in $L^2(\R^n)$ (and hence in $L^2_r(\R^n)$);
by the above the sequence $(A\phi_k)$ is Cauchy in $L^2(\R^n)$, and
$A\phi_k \to A\phi$ in $L^2_r(\R^n)$; this implies that $A\phi \in
L^2(\R^n)$ and $\|A\phi\|/\|\phi\| \leq \|A\|$. We now conclude that
$\hat{\sigma}(A) = A|_{L^2(\R^n)}$ is a bounded operator on $L^2(\R^n)$ and
$\|\hat{\sigma}(A)\| \leq \|A\|$. Injectivity of $\hat{\sigma}$ follows from
the fact that $L^2(\R^n)$ is dense in $L^2_r(\R^n)$, and $\hat{\sigma}$ is
trivially an algebra homomorphism. Next, it is clear that
$\hat{\sigma}(T_{\vec{a}}) = \widehat{M}_{e^{\iota\vec{a}\cdot\vec{x}}}
\in \widehat{L^\infty(\R^n)}$ for any $\vec{a} \in \R^n$. But any
$A \in L^\infty_r(\R^n)$ commutes with all $T_{\vec{a}}$, and
hence $\hat{\sigma}(A)$ commutes with all $\hat{\sigma}(T_{\vec{a}})$, so
we conclude $A \in \widehat{L^\infty(\R^n)}$ since the latter is a maximal
abelian subalgebra of $B(L^2(\R^n))$ and is generated as a von Neumann
algebra by the operators $\hat{\sigma}(T_{\vec{a}})$. Weak*-continuity
of $\hat{\sigma}$ holds for bounded nets since the weak* topology agrees
with the weak operator topology on bounded sets, and it then holds in
general by the Krein-Smulian theorem.\hfill\hal
\bigskip

\noindent {\bf Definition 17.} Define a map $\sigma: L^\infty_r(\R^n)
\to L^\infty(\R^n)$ by composing $\hat{\sigma}$ with the natural isomorphism
of $L^\infty(\R^n)$ with $\widehat{L^\infty(\R^n)}$. We call $\sigma$ the
{\it symbol map}.
\bigskip

We will retain both notations $\sigma$ and $\hat{\sigma}$, as
the latter will occasionally continue to be useful.
\medskip

We give three easy corollaries of Theorem 16.
\bigskip

\noindent {\bf Corollary 18.} {\it $L^\infty_r(\R^n)$ is a maximal
abelian subalgebra of $B(L^2_r(\R^n))$.}
\bigskip

(It is abelian since $\phi$ is injective, and maximality is then an
immediate consequence of Definition 13 (a).)
\bigskip

\noindent {\bf Corollary 19.} {\it Let $S \subset \R^n$ be measurable
and let $H_S \subset L^2_r(\R^n)$ be as in Definition 4. Then
$H_S$ is invariant for every $A \in L^\infty_r(\R^n)$, and the
distributional equation
$$A\phi = (\sigma(A)\check{\phi})\hat{\phantom{i}}$$
holds for all $\phi \in H_S$.}
\bigskip

This holds since the restriction of $A$ to $L^2(\R^n)$ is multiplication
by $\sigma(A)$ in the untransformed picture. It follows that
$\widehat{L^2(S)}$ is invariant for $A|_{L^2(\R^n)}$, and hence $H_S$
is invariant for $A$ by continuity. The equation for $A\phi$ also holds
by continuity.
\medskip

An algebra is {\it ergodic} for a group action if the only fixed points
of the action are scalar multiples of the identity. As the action of
$\R^n$ on $L^\infty(\R^n)$ by translations is ergodic and pulls
back under $\sigma$ (which is injective) to the action $\gamma$ on
$L^\infty_r(\R^n)$, the following corollary is immediate.
\bigskip

\noindent {\bf Corollary 20.} {\it The action of $\gamma$ on
$L^\infty_r(\R^n)$ is ergodic.}
\bigskip

The preceding result is in contrast
with the von Neumann algebra discussed in the last section, which
contains all bounded multiplication operators on $L^2(\R^n)$,
and hence is vastly non-ergodic for this action.
\medskip

Next we give a more serious corollary.
\bigskip

\noindent {\bf Corollary 21.} {\it The sets $\{K_f: f \in \S(\R^n)\}$
and ${\rm span}\{T_{\vec{a}}: \vec{a} \in \R^n\}$ are strong operator
dense in $L^\infty_r(\R^n)$.}
\medskip

\noindent {\it Proof.} It will suffice to verify that the $K_f$
are dense in $L^\infty_r(\R^n)$.
\medskip

Let $f$ be a smooth, compactly supported function on $\R^n$ such that
$\int f = 1$. Then define a sequence $(f_k)$ by $f_k(\vec{p})
= k^n f(k{\vec{p}})$.
\medskip

Let $A \in L^\infty_r(\R^n)$. We will show that $A$ is strong operator
approximated by convolution operators by Schwartz functions. Observe first
that the sequence $(K_{f_k})$ strong operator converges to the identity
operator on $L^2_r(\R^n)$, since $\|K_{f_k}\| \to 1$ (a consequence of
the elementary estimate [5, Theorem II.1.6]; see Lemma 36) and
$K_{f_k}g \to g$ in $\S(\R^n)$ for any
Schwartz function $g$. So $K_{f_k}A \to A$ strong operator and it will
suffice to approximate $K_{f_k}A$ for arbitrary $k$. Thus, fix $k$ and let
$$A_j = \int_{\R^n} f_j(\vec{p})\gamma_{\vec{p}}(K_{f_k}A)\, d\vec{p}.$$
This integral can be defined weakly, i.e., by the formula
$$\langle A_j \phi, \psi\rangle = \int_{\R^n} f_j(\vec{p})
\langle \gamma_{\vec{p}}(K_{f_k}A)\phi, \psi\rangle\, d\vec{p}$$
for $\phi, \psi \in L^2_r(\R^n)$. Since $\gamma$ is implemented by unitary
operators on $L^2_r(\R^n)$, it follows that $\|A_j\| \leq \|f\|_1\|K_{f_k}A\|$.
Now if $\sigma(A) = g \in L^\infty(\R^n)$, then a classical computation shows
that the restriction of $K_{f_k}A$ to $L^2(\R^n) \subset L^2_r(\R^n)$ is
$\widehat{M}_{\check{f}_k g}$, whereas the restriction of $A_j$ is
$\widehat{M}_{f_j * (\check{f}_k g)}$. It therefore follows that
$A_j \to K_{f_k}A$ strongly on vectors in $L^2(\R^n)$, and since the latter
is a dense subspace of $L^2_r(\R^n)$ and the sequence $(A_j)$ is bounded,
we conclude that $A_j \to K_{f_k}A$ strongly. Moreover,
$\hat{f}_j(f_k * \hat{g}) \in \S(\R^n)$
and $A_j = K_{\hat{f}_j(f_k * \hat{g})}$ (on $L^2(\R^n)$, and hence by
continuity on $L^2_r(\R^n)$), so we conclude that $A$ is a strong
operator limit of convolution operators by Schwartz functions.\hfill\hal
\bigskip

Next, it is interesting to note that although the algebra $L^\infty_r(\R^n)$
is not closed under adjoints, its image under $\sigma$ is. Thus,
$L^\infty_r(\R^n)$ is a Banach $*$-algebra in a natural way.
\bigskip

\noindent {\bf Proposition 22.} {\it The subalgebra
$\sigma(L^\infty_r(\R^n)) \subset L^\infty(\R^n)$ is self-adjoint.}
\medskip

\noindent {\it Proof.} Let $F: L^2_r(\R^n) \to L^2_r(\R^n)$ be the
antiunitary flip map $F\phi(\vec{p}) = \overline{\phi(-\vec{p})}$ and
observe that $F\lambda T_{\vec{a}}F = \bar{\lambda}T_{-\vec{a}}$ for
any $\lambda \in \C$ and $\vec{a} \in \R^n$. It follows that
$L^\infty_r(\R^n)$ is stable under conjugation by $F$.
Regarding $F$ also as an operator on $L^2(\R^n) \subset L^2_r(\R^n)$, we
have $FBF = B^*$ for any $B \in \widehat{L^\infty(\R^n)} \subset
B(L^2(\R^n))$. We conclude that
$$\hat{\sigma}(A)^* = F\hat{\sigma}(A)F = \hat{\sigma}(FAF)
\in \hat{\sigma}(L^\infty_r(\R^n))$$
for any $A \in L^\infty_r(\R^n)$.\hfill\hal
\bigskip

The preceding proof shows that the pullback of the involution on
$L^\infty(\R^n)$ to $L^\infty_r(\R^n)$ is given by conjugation with
$F$. Thus, it is isometric on $L^\infty_r(\R^n)$.
\medskip

We close this section with a problem on the existence of idempotents
in $L^\infty_r(\R^n)$.
\bigskip

\noindent {\bf Remarks 23.}
\smallskip

\noindent {\bf (a)} The map $\sigma: L^\infty_r(\R^n) \to L^\infty(\R^n)$
is not surjective. If it were, then by the open mapping theorem it would
be an isomorphism and the inverse map would be bounded. However,
$\|\sigma(T_{\vec{a}})\| = 1$ for all $\vec{a} \in \R^n$, while
$\|T_{\vec{a}}\| \to \infty$ as $|\vec{a}| \to \infty$, so the inverse
map is unbounded.
\medskip

\noindent {\bf (b)} It follows that $\sigma(L^\infty_r(\R^n))$ does
not contain every projection in $L^\infty(\R^n)$, since these span
$L^\infty(\R^n)$. In fact, one can show by direct computation that
in the case $n = 1$ the characteristic function of any interval
(besides the empty interval and the entire real line) fails to belong
to $\sigma(L^\infty_r(\R))$.
\medskip

We sketch the construction. First, if the characteristic function of
any half-infinite interval belonged to $\sigma(L^\infty_r(\R))$ then
so would any translation of it, and hence, taking differences, the
characteristic function of some bounded interval would also belong
to $\sigma(L^\infty_r(\R))$. Thus, it suffices to consider the bounded
case. Now the Fourier transform $f$ of any
bounded interval $J$ decays like $1/p$. We must find $\phi \in L^2_r(\R)$
such that $\|f * \phi\|_r/\|\phi\|_r$ is arbitrarily large. For
$N$ large and $a, b > 0$ depending on $f$ we let $\phi$ be the
characteristic function of a disjoint union of $N$ intervals $I_j$
$(1 \leq j \leq N)$, each of length $b$ and distance $ja$ from th
 origin, such that $f \approx c/ja$ around the interval $-I_j$ for some
nonzero complex scalar $c$, and hence $f * \chi_{I_j}$ is approximately
$bc/ja$ on $[0,b]$. (Here $\chi$ denotes characteristic function.)
Then $\|\phi\|_r^2$ is roughly $(b/a)\sum_1^N 1/j$ and $\|f * \phi\|_r^2$
is roughly at least $(|c|^2b^3/a^2)(\sum_1^N 1/j)^2$.
As $N$ goes to infinity, the divergence of $\sum_1^\infty 1/j$ implies
that $\|f * \phi\|_r/\|\phi\|_r$ goes to infinity. So convolution by $f$
is not bounded on $L^2_r(\R)$.
\bigskip

\noindent {\bf Problem 24.} Does $L^\infty_r(\R^n)$ contain any
nontrivial idempotents?
\bigskip
\bigskip

\noindent {\bf 5. Other non self-adjoint algebras}
\bigskip

Using the basic ingredients introduced in previous sections (the algebra
$L^\infty_r(\R^n)$, the action $\gamma$, the operators $T_{\vec{a}}$
and $K_f$) one can identify a large variety of operator algebras
which are analogous to various classical function algebras. We do
this now. In each case, it is easy to see (or vacuous) that the
analogous classical construction produces the analogous function
algebra.
\bigskip

\noindent {\bf Definition 25.} Let $T_{\vec{a}}$ and $K_f$ be as in
Definition 1 and $L^\infty_r(\R^n)$ and $\gamma$ as in Definition 13.
In all cases except (e) and (f), we give the following spaces the
operator norm they inherit from $L^\infty_r(\R^n) \subset B(L^2_r(\R^n))$.
\medskip

\noindent {\bf (a)} Let $C_{0,r}(\R^n)$ be the norm closure of the set
$\{K_f: f \in \S(\R^n)\}$.
\medskip

\noindent {\bf (b)} Let $AP_r(\R^n)$ be the norm closure of
${\rm span}\{T_{\vec{a}}: \vec{a} \in \R^n\}$.
\medskip

\noindent {\bf (c)} Let $UC_r(\R^n)$ be the set of $A \in L^\infty_r(\R^n)$
such that the function $\vec{a} \mapsto \gamma_{\vec{a}}(A)$ is continuous
for the norm topology.
\medskip

\noindent {\bf (d)} Let $C_{b,r}(\R^n)$ be the set of operators
$A \in L^\infty_r(\R^n)$ such that $B \in C_{0,r}(\R^n)$ implies
$AB \in C_{0,r}(\R^n)$.
\medskip

\noindent {\bf (e)} Define $C^k_{0,r}(\R^n)$ inductively by setting
$C^0_{0,r}(\R^n) = C_{0,r}(\R^n)$ and letting $C^{k+1}_{0,r}(\R^n)$
be the set of operators $A \in C_{0,r}(\R^n)$ such that the norm limit
$$D_i A = \lim_{t \to 0} {{1}\over{t}}(\gamma_{t\vec{e}_i}(A) - A)$$
($1 \leq i \leq n$) exists and belongs to $C^k_{0,r}(\R^n)$.
We give $C^{k+1}_{0,r}(\R^n)$ the norm
$$\|A\|_{k+1} = \max(\|A\|, \|D_1A\|_k, \ldots, \|D_nA\|_k).$$

\noindent {\bf (f)} Let ${\rm Lip}_r(\R^n)$ be the set of operators
$A \in L^\infty_r(\R^n)$ such that
$$L(A) = \sup_{\vec{a} \in \R^n} {{1}\over{|\vec{a}|}}
\|\gamma_{\vec{a}}(A) - A\|$$
is finite (using the convention $0/0 = 0$). We give it the norm
$\|A\|_L = \max(\|A\|, L(A))$.
\medskip

\noindent {\bf (g)} Let $H^\infty_r(\R)$ be the strong operator closure
of ${\rm span}\{T_a: a \geq 0\} \subset L^\infty_r(\R)$.
\bigskip

This list could easily be extended further.
\medskip

$AP$ stands for ``almost periodic,'' $UC$ for ``uniformly continuous,''
$C_b$ for ``continuous and bounded,'' and ${\rm Lip}$ for ``Lipschitz.''
\medskip

The next proposition is straightforward and we omit its proof.
\bigskip

\noindent {\bf Proposition 26.} {\it Each of the spaces in Definition 25
is a Banach algebra. The images of all but $H^\infty_r(\R^n)$ under
$\sigma$ are self-adjoint subalgebras of $L^\infty(\R^n)$.}
\bigskip

(Note that in cases (e) and (f) the Banach algebra law holds only
in its weak form, i.e., $\|xy\| \leq C\|x\|\|y\|$ for some constant
$C$. Of course, this occurs classically as well.)
\medskip

Next we observe that the expected containments hold. This uses a lemma
which is of independent interest.
\bigskip

\noindent {\bf Lemma 27.} {\it Let $f \in \S(\R^n)$
be compactly supported and satisfy $\int f = 1$.
Then $A \in L^\infty_r(\R^n)$ belongs to $C_{0,r}(\R^n)$ if and only if
$A \in UC_r(\R^n)$ and $K_{f_k}A \to A$ in norm, where $f_k(\vec{p}) = 
k^n f(k{\vec{p}})$.}
\medskip

\noindent {\it Proof.} Observe first that $K_{f_k}K_g = K_{f_k * g} \to K_g$
in norm for all $g \in \S(\R^n)$; this follows from the fact that
$f_k * g \to g$ in $\S(\R^n)$. Moreover, using the fact that $\|K_{f_k}\|
\to 1$, an $\epsilon/3$ argument shows that the set of $A$ which satisfy
$K_{f_k}A \to A$ is closed in norm. This proves that every $A \in
C_{0,r}(\R^n)$ satisfies $K_{f_k}A \to A$. Also, direct computation shows
that $\gamma_{\vec{a}}(K_f) = K_{e^{-\iota\vec{a}\cdot{\vec{p}}}f}$, and
this is norm continuous in $\vec{a}$ for any $f \in \S(\R^n)$ by
[5, Theorem II.1.6] (cf.\ Lemma 36). Using
the fact that $UC_r(\R^n)$ is norm closed, it follows that $C_{0,r}(\R^n)$
is contained in $UC_r(\R^n)$. This completes the proof of the forward
implication.
\medskip

For the reverse implication, let $A \in UC_r(\R^n)$ and suppose
$K_{f_k}A \to A$ in norm. It will suffice to show that $K_{f_k}A
\in C_{0,r}(\R^n)$ for arbitrary $k$. Thus, fix $k$ and let
$$A_j = \int_{\R^n} f_j(\vec{p})\gamma_{\vec{p}}(K_{f_k}A)\, d\vec{p}$$
as in the proof of Corollary 21. Now $K_{f_k}A =
\int f_k(\vec{p})K_{f_k}A\, d\vec{p}$ so
$$\|A_j - K_{f_k}A\| \leq \int_{\R^n} |f_j(\vec{p})|
\|\gamma_{\vec{p}}(K_{f_k}A) - K_{f_k}A\|\, d\vec{p}.$$
This goes to zero as $j \to \infty$ since $K_{f_k}A \in UC_r(\R^n)$ implies
that the function $\vec{p} \mapsto \gamma_{\vec{p}}(K_{f_k}A) - K_{f_k}A$
is continuous in norm, and it is clearly bounded
and zero at $\vec{p} = 0$. However, $A_j = K_{\hat{f}_j(f_k * \hat{g})}$
where $g = \sigma(A) \in L^\infty(\R^n)$; and
$\hat{f}_j(f_k * \hat{g}) \in \S(\R^n)$, so $A_j$ is of the form
$K_f$ for $f \in \S(\R^n)$. We conclude that $K_{f_k}A$, and hence
$A$, belongs to $C_{0,r}(\R^n)$.\hfill\hal
\bigskip

The condition $A \in UC_r(\R^n)$ in Lemma 27 can be weakened to
$A \in C_{b,r}(\R^n)$. In fact the proof of the backward implication
becomes easier, for then we know immediately that $K_{f_k}A \in
C_{0,r}(\R^n)$. The point is that proving $UC_r(\R^n) \subset
C_{b,r}(\R^n)$ is tantamount to showing that $A \in UC_r(\R^n)$
implies $K_{f_k}A \in C_{0,r}(\R^n)$.
\bigskip

\noindent {\bf Proposition 28.} {\it We have
$$AP_r(\R^n), {\rm Lip}_r(\R^n), C_{0,r}(\R^n)
\subset UC_r(\R^n) \subset C_{b,r}(\R^n).$$}

\noindent {\it Proof.} We showed that $C_{0,r}(\R^n) \subset UC_r(\R^n)$
in the lemma. The proof of $AP_r(\R^n) \subset UC_r(\R^n)$ is similar,
now using the fact that $\gamma_{\vec{a}}(T_{\vec{b}}) =
e^{-\iota\vec{a}\cdot\vec{b}} T_{\vec{b}}$ is norm continuous in $\vec{a}$.
The containment ${\rm Lip}_r(\R^n) \subset UC_r(\R^n)$ holds because
$$\|\gamma_{\vec{a}}(A) - \gamma_{\vec{b}}(A)\| =
\|\gamma_{\vec{b}}(\gamma_{\vec{a} - \vec{b}}(A) - A)\|
\leq |\vec{a} - \vec{b}|L(A)$$
for all $\vec{a}, \vec{b} \in \R^n$ and all $A \in {\rm Lip}_r(\R^n)$.
\medskip

For the final containment, let $A \in UC_r(\R^n)$ and
$B \in C_{0,r}(\R^n)$; we must show that $AB \in C_{0,r}(\R^n)$. As
above, $B \in UC_r(\R^n)$, so $AB \in UC_r(\R^n)$ by Proposition
26. Also $K_{f_k}B \to B$ implies $K_{f_k}AB = AK_{f_k}B \to AB$.
So $AB \in C_{0,r}(\R^n)$ by the lemma.\hfill\hal
\bigskip

The next result is also basic.
\bigskip

\noindent {\bf Proposition 29.} {\it ${\rm Lip}_r(\R^n)$ and
$H^\infty_r(\R)$ are dual spaces.}
\medskip

\noindent {\it Proof.} It is standard that any strong operator closed
subspace of $B(H)$ is weak operator closed, and hence weak* closed, as
in the proof of Proposition 15. This shows that $H^\infty_r(\R)$ is
a dual space. For ${\rm Lip}_r(\R^n)$, let
$\M = \bigoplus_{\vec{a}} L^\infty_r(\R^n)$ be the $l^\infty$ direct sum
over $\vec{a} \in \R^n$, $\vec{a} \neq 0$, and consider the map
$d: {\rm Lip}_r(\R^n) \to \M$ defined by $(dA)_{\vec{a}} =
(\gamma_{\vec{a}}(A) - A)/|\vec{a}|$. Then $\|A\|_L =
{\rm max}(\|A\|, \|dA\|)$, i.e., $\|\cdot\|_L$ is the graph norm.
As $\gamma_{\vec{a}}$ is weak* continuous, it easily follows
that the graph of $d$ is weak* closed in $L^\infty_r(\R^n)
\oplus \M$, so this graph, which is isometric to ${\rm Lip}_r(\R^n)$,
is a dual space.\hfill\hal
\bigskip

We note that according to [4, Proposition 3.1.23], ${\rm Lip}_r(\R^n)$
is equivalently the set of operators $A$ in $L^\infty_r(\R^n)$ such that
the partial derivatives $D_iA$ ($1 \leq i \leq n$) as in Definition 25
(e) exist as weak operator limits.
\medskip

We now turn to the relation of the above spaces to the symbol map $\sigma$.
The first result is easy.
\bigskip

\noindent {\bf Proposition 30.} {\it Let $X_r$ be any of the spaces
in Definition 25 and let $X$ be its classical analog. Then
$\sigma(X_r) \subset X$.}
\bigskip

The unboundedness argument of Remark 23 (a) applies to every algebra
in Definition 25 to show that $\sigma(X_r)$ must be a proper subset of
$X$. A basic question is whether
$\sigma(X_r) = X \cap \sigma(L^\infty_r(\R^n))$. We first show that this
holds for $H^\infty_r(\R)$.
\bigskip

\noindent {\bf Proposition 31.} {\it $H^\infty_r(\R)$ consists of
precisely those operators $A \in L^\infty_r(\R)$ such that
$$\inf {\rm supp}\, A\phi \geq \inf {\rm supp}\, \phi$$
for all $\phi \in L^2_r(\R)$ with support bounded from below.
We have
$$\sigma(H^\infty_r(\R)) = H^\infty(\R) \cap \sigma(L^\infty_r(\R)).$$}

\noindent {\it Proof.} It is clear that every $T_a$, $a \geq 0$, translates
supports of functions in $L^2_r(\R)$ to the right. It follows that any
finite linear combination of the $T_{\vec{a}}$ respects lower bounds of
supports, and, taking strong operator limits, that this is true of
any operator in $H^\infty(\R)$. Conversely, let $A \in L^\infty_r(\R)$
and suppose $\inf {\rm supp}\, A\phi \geq \inf {\rm supp}\, \phi$ for all
$\phi \in L^2_r(\R)$ with support bounded from below. Then this is
true in particular for all $\phi \in L^2(\R) \subset L^2_r(\R)$,
so that $\hat{\sigma}(A) \in \widehat{H^\infty(\R)} \subset
\widehat{L^\infty(\R)}$. Adopting the notation of the proof of
Corollary 21, we have that $A$ is strong operator approximated by
operators of the form $K_{\hat{f}_j(f_k * \hat{g})}$ where $\sigma(A) =
g \in H^\infty(\R)$. Since the $f_k$ have compact support,
the support of $f_k$ must converge to zero as $k \to \infty$,
and it follows that $A$ is strong operator approximated
by operators of the form $K_h$ with ${\rm supp}\, h \subset [0,\infty)$.
This implies that $A \in H^\infty_r(\R)$.
\medskip

The equation $\sigma(H^\infty_r(\R)) =
H^\infty(\R) \cap \sigma(L^\infty_r(\R))$ follows easily.\hfill\hal
\bigskip

Next, we have partial positive results to the same question for
$C_{0,r}(\R^n)$ and $AP_r(\R^n)$.
\bigskip

\noindent {\bf Lemma 32.} {\it If $f \in L^\infty(\R^n)$ and
$\hat{f}$ has compact support then $f \in \sigma(L^\infty_r(\R^n))$.
Moreover, if $S \subset \R^n$ is compact then there is a constant
$C > 0$ such that $\|\sigma^{-1}(f)\| \leq C\cdot \|f\|_\infty$
for all $f \in L^\infty(\R^n)$ with ${\rm supp}\, \hat{f} \subset S$.}
\medskip

\noindent {\it Proof.} Fix a compact set $S \subset \R^n$ and let
$\CC_1 = [-N/2, N/2]^n$ be a cube which contains $S$. Then the cubes
$\CC_1 + N\cdot \Z^n$ tile $\R^n$. Let $\CC_2, \ldots, \CC_{3^n}$
be the cubes adjacent to $\CC_1$ in this tiling.
\medskip

Fix an arbitrary $f \in L^\infty(\R^n)$ such that ${\rm supp}\, \hat{f}
\subset S$ (in the distributional sense). Regard $\widehat{M}_f$ as
acting either on $L^2(\R^n)$ or, perhaps unboundedly, on $L^2_r(\R^n)$.
For a cube $\CC \subset \R^n$ of any size and location, let $P_\CC:
\phi \mapsto \chi_\CC\cdot \phi$ be the orthogonal projection given by
restriction to $\CC$ ($\chi$ denotes characteristic function)
acting on either $L^2(\R^n)$ or $L^2_r(\R^n)$. We claim that there
is a cube $\CC$ whose side has length $3N$ such that the norm of
$P_\CC \widehat{M}_f P_\CC$ acting on $L^2_r(\R^n)$ is at least
$3^{-n}$ times the norm of $\widehat{M}_f$ acting on $L^2_r(\R^n)$
(or infinite if the latter is infinite). To see this, for
$1 \leq i \leq 3^n$ let $A_i = \sum P_\CC \widehat{M}_f P_\CC$, taking
the sum over all cubes $\CC$ belonging to the tiling
$\CC_i' + 3N\cdot \Z^n$ where $\CC_i'$ is the cube with side length
$3N$ consisting of $\CC_i$ and its neighbors.
Now ${\rm supp}\, \hat{f} \subset \CC_1$
implies that for any $\phi \in L^2_r(\R^n)$ supported outside $\CC_i'$
we must have $\widehat{M}_f \phi = \hat{f}*\phi = 0$ on $\CC_i$. Thus,
$$P_{\CC_i}\widehat{M}_f P_{\CC_i'}\phi = P_{\CC_i}\widehat{M}_f \phi,$$
that is, $A_i\phi = \widehat{M}_f \phi$ on $\CC_i$,
for any $\phi \in L^2_r(\R^n)$; more generally, $A_i\phi$ must agree with
$\widehat{M}_f\phi$ on the set $\CC_i + 3N\cdot\Z^n$. Since these
$3^n$ sets ($1 \leq i \leq 3^n$) tile $\R^n$, it follows that
$$\|A_1\|^2 + \cdots + \|A_{3^n}\|^2 \geq \|\widehat{M}_f\|^2$$
on $L^2_r(\R^n)$, and hence $\|A_i\| \geq 3^{-n/2}\|\widehat{M}_f\|$
(or is infinite, if $\|\widehat{M}_f\|$ is infinite)
for some $i$. But $\|A_i\|$ is the supremum of the norms
$\|P_\CC \widehat{M}_f P_\CC\|$ where $\CC$ ranges over the sum which
defines $A_i$. Thus, there is some cube $\CC$ whose side has
length $3N$ such that $\|P_\CC \widehat{M}_f P_\CC\| \geq
3^{-n}\|\widehat{M}_f\|$, as claimed.
\medskip

Find a constant $C > 0$ such that if $\CC$ is any cube whose side
has length $3N$, and $a$ and $b$ are respectively the maximum and
minimum of the function $E(\vec{p})$ on $\CC$, then $a/b \leq C$.
This is possible because the ratio goes to 1 as the cube goes to
infinity. Now, continuing to let $f \in L^\infty(\R^n)$ be arbitrary
such that ${\rm supp}\, \hat{f} \subset S$, find a cube $\CC$ so that
$\|P_\CC \widehat{M}_f P_\CC\| \geq 3^{-n}\|\widehat{M}_f\|$ on
$L^2_r(\R^n)$, and let $\phi \in L^2_r(\R^n)$ be supported on $\CC$.
Let $a$ and
$b$ be the respective maximum and minimum of $E(\vec{p})$ on this
cube. Then the norm of $\phi$ in $L^2(\R^n)$ is related to its
norm in $L^2_r(\R^n)$ by $\|\phi\|_r \geq \sqrt{b}\|\phi\|$.
Similarly, we have $\|P_\CC\widehat{M}_f \phi\|_r \leq
\sqrt{a}\|P_\CC\widehat{M}_f \phi\|$. But $\|P_\CC\widehat{M}_f \phi\|
\leq \|f\|_\infty \|\phi\|$, so that
$$\|P_\CC\widehat{M}_f \phi\|_r \leq \sqrt{a}\|P_\CC\widehat{M}_f \phi\|
\leq \sqrt{a}\|f\|_\infty \|\phi\|
\leq \sqrt{C}\|f\|_\infty\|\phi\|_r.$$
Together with the choice of $\CC$, we conclude that the norm of
$\widehat{M}_f = \sigma^{-1}(f)$ on $L^2_r(\R^n)$ is at most
$3^n\sqrt{C}\|f\|_\infty$. This completes the proof.\hfill\hal
\bigskip

Note that if $\hat{f}$ has compact support then $f$ is infinitely
differentiable [17, Theorem 7.23]. This aspect of the above result
will be strengthened in Corollary 37.
\bigskip

\noindent {\bf Corollary 33.} {\it Let $X_r$ be any of the spaces
of Definition 25 and let $X$ be its classical counterpart. Then
$\sigma(X_r)$ contains every $f \in X$ such that $\hat{f}$ has
compact support.}
\bigskip

\noindent {\bf Theorem 34.} {\it We have
$$\sigma(C_{0,r}(\R^n)) = C_0(\R^n) \cap \sigma(UC_r(\R^n))$$
and
$$\sigma(AP_r(\R^n)) = AP(\R^n) \cap \sigma(UC_r(\R^n)).$$}

\noindent {\it Proof.} In both cases, one containment follows
from Propositions 28 and 30. For the other direction, choose $g$ in
$C_0(\R^n) \cap \sigma(UC_r(\R^n))$ (respectively,
$AP(\R^n) \cap \sigma(UC_r(\R^n))$) and let $A = \sigma^{-1}(g)$. Choose
$f \in \S(\R^n)$ such that $\int f = 1$ (i.e., $\hat{f}(0) = 1$) and
$\hat{f}$ has compact support, and define $f_k(\vec{p}) =  k^n f(k{\vec{p}})$.
Then $A \in UC_r(\R^n)$ implies that the sequence
$$A_j = \int_{\R^n} f_j(\vec{p})\gamma_{\vec{p}}(A)\, d\vec{p}$$
converges to $A$ in norm, as in the proof of Lemma 27. Moreover,
$A_j = K_{\hat{f}_j\hat{g}}$, so that $\sigma(A_j)\hat{\phantom{i}}
= \hat{f}_j\hat{g}$ has compact support, while $f_j * g \in C_0(\R^n)$
(respectively $AP(\R^n)$) is clear. Corollary 33 then implies
that $A_j$ is in $C_{0,r}(\R^n)$ (respectively, $AP_r(\R^n)$).
Taking $j \to \infty$, we find that $A$ belongs to the same space,
as desired.\hfill\hal
\bigskip

\noindent {\bf Corollary 35.} {\it $C_{0,r}(\R^n)^+$ and $AP_r(\R^n)$
are inverse closed subalgebras of $UC_r(\R^n)$.}
\bigskip

We now proceed to a negative result which states that there exists a
function in $\sigma(C_{0,r}(\R^n))$ which belongs to
$C_0^d(\R^n)$ but not $\sigma(C_{0,r}^d(\R^n))$.
\bigskip

\noindent {\bf Lemma 36.} {\it Let $f: \R^n \to \C$ be measurable
and suppose
$$\int_{\R^n} (1 + |\vec{p}|^{1/2})|f(\vec{p})|\, d\vec{p} < \infty.$$
Then $K_f$, the operator of convolution by $f$ on $L^2_r(\R^n)$, belongs
to $L^\infty_r(\R^n)$ and
$$\|K_f\| \leq \int_{\R^n}
\left(1 + {{|\vec{p}|}\over{m}} + {{|\vec{p}|^2}\over{m^2}}\right)^{1/4}
|f(\vec{p})|\, d\vec{p}.$$}

\noindent {\it Proof.} It is enough to show that $K_f$ is bounded on
$L^2_r(\R^n)$; commutation with the operators $T_{\vec{a}}$ is then
trivial and implies $K_f \in L^\infty_r(\R^n)$. 
\medskip

For the sharpest bound on $\|K_f\|$, we pass to the operator
$\tilde{K}_f = M_{\sqrt{E}}^{-1}K_fM_{\sqrt{E}}$ on $L^2(\R^n)$ as in
Section 3. Recall that
$$\tilde{K}_f \phi(\vec{p}) =
\int f(\vec{p}\,{}') \sqrt{{E(\vec{p} - \vec{p}\,{}')}\over{E(\vec{p})}}
\phi(\vec{p} - \vec{p}\,{}')\, d\vec{p}\,{}' = 
\int f(\vec{p} - \vec{p}\,{}') \sqrt{{E(\vec{p}\,{}')}\over{E(\vec{p})}}
\phi(\vec{p}\,{}')\, d\vec{p}\,{}',$$
so $\tilde{K}_f$ is an integral operator with kernel
$$K(\vec{p}, \vec{p}\,{}') =
f(\vec{p} - \vec{p}\,{}') \sqrt{E(\vec{p}\,{}')/E(\vec{p})}.$$
We claim that the integrals $\int |K(\vec{p}, \vec{p}\,{}')|\, d\vec{p}$ and
$\int |K(\vec{p}, \vec{p}\,{}')|\, d\vec{p}\,{}'$ are uniformly bounded by
$$\int_{\R^n}
\left(1 + {{|\vec{p}|}\over{m}} + {{|\vec{p}|^2}\over{m^2}}\right)^{1/4}
|f(\vec{p})|\, d\vec{p}.$$
By [5, Theorem II.1.6] the norm of $\widetilde{K}_f$, and hence of
$K_f$, is then less than or equal to this uniform bound. The claim
follows by a short computation using the elementary inequality
$$\sqrt{{E(\vec{p} - \vec{p}\,{}')}\over{E(\vec{p}\,{}')}} \leq
\left(1 + {{|\vec{p}|}\over{m}} + {{|\vec{p}|^2}\over{m^2}}\right)^{1/4}$$
and the assumption on integrability of $f$ then implies that the bound
on $K_f$ is finite.\hfill\hal
\bigskip

The argument in the above proof, while easier than the one in
the proof of Lemma 32, does not apply there because in that
situation $\hat{f}$ (which plays the role played by $f$ here)
is in general not a function.
\bigskip

\noindent {\bf Corollary 37.} {\it Let $f$ be a compactly supported,
$d$-times continuously differentiable function on $\R^n$, where
$d > (n+1)/2$. Then $f \in \sigma(C_{0,r}(\R^n))$.}
\medskip

\noindent {\it Proof.} Observe first that
$$\int |\vec{p}|^{1/2}|\hat{f}(\vec{p})|\, d\vec{p}
\leq \left(\int {{d\vec{p}}\over{1 + |\vec{p}|^{2d-1}}}\right)^{1/2}
\left(\int (1 + |\vec{p}|^{2d - 1})|\vec{p}| |\hat{f}(\vec{p})|^2\, d\vec{p}\right)^{1/2}$$
by the Cauchy-Schwarz inequality. The first integral on the right
is finite by our choice of $d$. For the second, observe that
$$\int |\vec{p}|^{2d}|\hat{f}(\vec{p})|^2\, d\vec{p}
\leq \int n^d \max_{1 \leq i \leq n} |p_i^d \hat{f}(\vec{p})|^2\, d\vec{p}
\leq n^d \sum_{1 \leq i \leq n}
\int\left|{{\partial^d f}\over{\partial x_i^d}}(\vec{x})\right|^2\, d\vec{x}
< \infty.$$
Since $\hat{f}$
is continuous, and hence bounded in a neighborhood of the origin, this
implies that $\int |\vec{p}||\hat{f}(\vec{p})|^2\, d\vec{p}$ is also finite,
hence $\int |\vec{p}|^{1/2}|\hat{f}(\vec{p})|\, d\vec{p}$ is finite and
therefore so is $\int |\hat{f}(\vec{p})|\, d\vec{p}$. We conclude that
$\hat{f}$ satisfies the hypothesis of Lemma 36. Thus $K_{\hat{f}} \in
L^\infty_r(\R^n)$, and to $f \in \sigma(L^\infty_r(\R^n))$.
\medskip

To show that $f \in \sigma(C_{0,r}(\R^n))$, by Theorem 34 it is enough to
show that $f \in \sigma(UC_r(\R^n))$, i.e., that $K_{\hat{f}} \in UC_r(\R^n)$.
Now
$$\gamma_{\vec{a}}(K_{\hat{f}}) - K_{\hat{f}}
= K_{(e^{-\iota\vec{a}\cdot{\vec{p}}} - 1)\hat{f}}.$$ 
Substituting $(e^{-\iota\vec{a}\cdot{\vec{p}}} - 1)\hat{f}$ for $f$
in the hypothesis of Lemma 36, the dominated convergence theorem implies
that the bound on $\|K_{(e^{-\iota\vec{a}\cdot{\vec{p}}} - 1)\hat{f}}\|$
given there goes to zero as $\vec{a} \to 0$. So
$\gamma_{\vec{a}}(K_{\hat{f}}) \to K_{\hat{f}}$ as $\vec{a} \to 0$. This
shows that $K_{\hat{f}} \in UC_r(\R^n)$, as desired.\hfill\hal
\bigskip

\noindent {\bf Theorem 38.} {\it Let $d > (n+1)/2$. Then
$\sigma(C_{0,r}^d(\R^n))$ is properly contained in
$C_0^d(\R^n) \cap \sigma(C_{0,r}(\R^n))$.}
\medskip

\noindent {\it Proof.} By Corollary 37, $\sigma(C_{0,r}(\R^n))$
contains every compactly supported function in $C_0^d(\R^n)$. Thus,
we need only show that there exists such a function which is not
in $\sigma(C_{0,r}^d(\R^n))$. Let $S$ be the unit ball in $\R^n$
and let $X$ be the set of operators $A$ in $C_{0,r}^d(\R^n)$ such
that $\sigma(A)$ is supported in $S$. We claim that $X$ is a
closed subspace of $C_{0,r}^d(\R^n)$. To see this, let $(A_k)$ be
a Cauchy sequence in $X$, say $A_k \to A \in C_{0,r}^d(\R^n)$, and
let $f \in \S(\R^n)$ be any Schwartz function which is constantly 1 on $S$.
Then $K_{\hat{f}}A_k = A_k$ for all $k$, so $K_{\hat{f}}A = A$ by
continuity. It follows that $f\cdot\sigma(A) = \sigma(A)$, and we deduce
that $\sigma(A)$ is supported on $S$, i.e., $A \in X$. Thus $X$ is
closed. Suppose $\sigma$ maps $X$ onto the closed subspace $Y$ of
$C_0^d(\R^n)$ consisting of all the functions in $C_0^d(\R^n)$ which
are supported on $S$. Then the open mapping theorem implies that
$\sigma$ is an isomorphism from $X$ onto $Y$.
\medskip

Let $g \in Y$ be arbitrary but nonzero. Then the norm of $g_k =
k^{-d}e^{\iota kx_1}g$ in $Y$ is bounded as $k \to \infty$, but
$\|D^d_1 \sigma^{-1}(g_k)\| \sim \sqrt{k}$.
This shows that $\sigma^{-1}$ is not bounded from $Y$ to $X$,
contradicting the conclusion of the last paragraph. Thus $\sigma(X)$,
and hence $\sigma(C_{0,r}^d(\R^n))$, does not contain all functions in
$C_0^d(\R^n)$ which are supported on $S$, as claimed.\hfill\hal
\bigskip

This still leaves our main question largely unanswered.
\bigskip

\noindent {\bf Problem 39.} For which relativistic spaces $X_r \subset
L^\infty_r(\R^n)$ with classical analog $X$ does $\sigma(X_r) =
X \cap \sigma(L^\infty_r(\R^n))$?
\bigskip
\bigskip

\noindent {\bf 6. Ideals and subalgebras}
\bigskip

We begin by identifying the maximal ideals of $C_{0,r}(\R^n)$.
\bigskip

\noindent {\bf Theorem 40.} {\it The maximal ideal space of
$C_{0,r}(\R^n)$ can be identified with $\R^n$, and the symbol
map $\sigma: C_{0,r}(\R^n) \to C_0(\R^n)$ can be identified
with the Gelfand transform.}
\medskip

\noindent {\it Proof.} For every $\vec{a} \in \R^n$, the map
$A \mapsto \sigma(A)(\vec{a})$ is a complex homomorphism on
$C_{0,r}(\R^n)$. Since $\sigma(C_{0,r}(\R^n))$ contains the
Schwarz functions on $\R^n$, it follows that the preceding
homomorphisms are distinct. We must show that every complex
homomorphism on $C_{0,r}(\R^n)$ is of this form.
\medskip

Let $\omega: C_{0,r}(\R^n) \to \C$ be a homomorphism, and let
$\A \subset C_{0,r}(\R^n)$ be the dense subalgebra consisting
of all operators of the form $K_f$ such that $f \in \S(\R^n)$
has compact support. We claim that $|\omega(K_f)| \leq \|f\|_1$
for every $K_f \in \A$. Letting $f^{*k}$ be the $k$th convolution
power of $f$, we have $K_f^k = K_{f^{*k}}$, and if $f$ is supported on
$[-N,N]^n$ then $f^{*k}$ is supported on $[-kN,kN]^n$, so as $k \to \infty$
$$\|K_f^k\| = \|f^{*k}\|_1\cdot O((nkN/m)^{1/2})
\leq \|f\|_1^k\cdot O(\sqrt{k})$$
by Lemma 36. We are using the fact that the dominant term in the bound on
$\|K_f^k\|$ given by Lemma 36 is $(|\vec{p}|^2/m^2)^{1/4} =
(|\vec{p}|/m)^{1/2}$. It follows that
$$|\omega(K_f)|^k = |\omega(K_f^k)| \leq C\sqrt{k}\|f\|_1^k,$$
and taking $k$th roots and letting $k \to \infty$ yields
$|\omega(K_f)| \leq \|f\|_1$, as claimed.
\medskip

If follows that $\omega$ extends by continuity to a complex homomorphism
on the convolution algebra $L^1(\R^n)$. These are well-known to be
given by point evaluations on $\R^n$ composed with the Fourier transform.
Thus $\omega$ applied to elements of $\A$ must be a point evaluation on
$C_0(\R^n)$ composed with $\sigma$, and hence this holds on $C_{0,r}(\R^n)$.
This completes the proof.\hfill\hal
\bigskip

Now we turn to more general ideals. We have the following examples
of ideals of $L^\infty_r(\R^n)$ and $C_{0,r}(\R^n)$.
\bigskip

\noindent {\bf Definition 41.} For any measurable $S \subset \R^n$,
let $\I(S)$ be the set of all $A \in L^\infty_r(\R^n)$
such that $\sigma(A)|_S \equiv 0$.
\bigskip

\noindent {\bf Proposition 42.} {\it Let $S \subset \R^n$ be measurable.
Then $\I(S)$ is a weak* closed ideal of $L^\infty_r(\R^n)$ and
$\I(S) \cap C_{0,r}(\R^n)$ is a norm closed ideal of $C_{0,r}(\R^n)$.}
\medskip

\noindent {\it Proof.} Let $T$ be the restriction map from $L^\infty(\R^n)$
to $L^\infty(\R^n - S)$. It is weak* continuous, hence $T \circ \sigma$
is weak* continuous, so $\I(S) = {\rm ker}(T \circ \sigma)$ is a weak*
closed ideal of $L^\infty_r(\R^n)$. The second assertion of the proposition
follows easily.\hfill\hal
\bigskip

\noindent {\bf Problem 43.} Are there any weak* closed ideals of
$L^\infty_r(\R^n)$ not of the form $\I(S)$, or any norm closed ideals
of $C_{0,r}(\R^n)$ not of the form $\I(S) \cap C_{0,r}(\R^n)$?
\bigskip

\noindent {\bf Remarks 44.}
\smallskip

\noindent {\bf (a)} If $\overline{S}$ is the closure of $S$, then
Proposition 30 implies that $\I(S) \cap C_{0,r}(\R^n) =
\I(\overline{S}) \cap C_{0,r}(\R^n)$. Conversely, if $S$ and $S'$ are
distinct closed sets then there is a Schwarz function which vanishes
on one but not the other, and hence $\I(S) \cap C_{0,r}(\R^n) \neq
\I(S') \cap C_{0,r}(\R^n)$.
\medskip

\noindent {\bf (b)} Taking advantage of Corollary 19, we can also
define spaces $L^\infty_r(S)$ and $C_{0,r}(S)$ as, respectively,
the strong operator and norm closures of the set
$\{K_f|_{H_S}: f \in \S(\R^n)\} \subset B(H_S)$.
\bigskip

Among the variety of possible constructions of other related operator algebras
(such as those just exhibited in Remark 44 (b)), perhaps the most interesting
are a family of algebras on the $n$-torus arising from the operators
$T_{\vec{a}}$ for $\vec{a} \in \Z^n$. Because harmonic analysis
techniques are generally more powerful on $\Z^n$ than on $\R^n$,
these algebras are somewhat more tractable than the algebras in
Sections 4 and 5: most of the results in those sections can be proven
more easily in the torus case. However, instead of giving direct proofs
we find it simpler to reduce to the real case.
\medskip

We adopt the convention that $\T^n = \R^n/2\pi\Z^n$.
\bigskip

\noindent {\bf Definition 45.}
\smallskip

\noindent {\bf (a)} Let $L^\infty_r(\T^n)$ be the strong operator closure
of ${\rm span}\{T_{\vec{a}}: \vec{a} \in \Z^n\}$ in $L^\infty_r(\R^n)$.
\medskip

\noindent {\bf (b)} Let $C_r(\T^n)$ be the norm closure of
${\rm span}\{T_{\vec{a}}: \vec{a} \in \Z^n\}$.
\medskip

\noindent {\bf (c)} Define $C^k_r(\T^n)$ inductively by setting
$C^0_r(\T^n) = C_r(\T^n)$ and letting $C^{k+1}_r(\T^n)$ be the
set of operators $A$ in $C^k_t(\T^n)$ such that $D_iA$ exists
as a norm limit, as in Definition 25 (e).
\medskip

\noindent {\bf (d)} Let ${\rm Lip}_r(\T^n)$ be the set of operators
$A$ in $L^\infty_r(\T^n)$ such that $L(A) < \infty$, as in
Definition 25 (f).
\medskip

\noindent {\bf (e)} Let $A_r(\T)$ be the norm closure of
${\rm span}\{T_a: a \in \N\}$.
\medskip

\noindent {\bf (f)} Let $H^\infty_r(\T)$ be the strong operator closure
of ${\rm span}\{T_a: a \in \N\}$.
\bigskip

\noindent {\bf Remarks 46.}
\smallskip

\noindent {\bf (a)} Composition with any invertible linear transformation
on $\R^n$ defines a bounded operator on $L^2_r(\R^n)$ (Remark 14 (c)).
Thus, if in the above definitions $\Z^n$ were replaced by any cocompact
lattice, the resulting spaces would have equivalent norms.
\medskip

\noindent {\bf (b)} Similarly, $L^2_r(\R^n)$ is naturally isomorphic
(but not isometric)
to the tensor product $l^2(\Z^n, 1/E(\vec{k})) \otimes L^2([0,1]^n)$.
The corresponding representation of the above algebras on this product
is trivial on the second factor, so it follows that they are all
isomorphically represented on $l^2(\Z^n, 1/E(\vec{k}))$ in the obvious way.
\bigskip

For $A \in L^\infty_r(\T^n)$, define the $N$th Ces\`aro mean of $A$ by
$$\tau_N(A) = {{1}\over{(2\pi)^n}}
\int_{\T^n} \gamma_{\vec{t}}(A) K_N^n(\vec{t}\,)\, d\vec{t},$$
where $K_N$ is the Fej\'er kernel
$$K_N(t) = \sum_{m = -N}^N\left(1 - {{|m|}\over{N+1}}\right)e^{\iota mt}
= {{1}\over{N+1}}\left({{\sin((N+1)t/2)}\over{\sin(t/2)}}\right)^2$$
($t \in \T$) and $K_N^n(\vec{t}\,) = K_N(t_1)\cdot\ldots\cdot K_N(t_n)$.
Observe that if $A$ is in ${\rm span}\{T_{\vec{a}}: \vec{a} \in \Z^n\}$
then $\tau_N(A)$ is a linear combination of the $T_{\vec{a}}$
with $\vec{a} \in \Z^n$ and $-N \leq a_i \leq N$ for all $i$.
By dominated convergence, the map $A \mapsto \tau_N(A)$ is weak*
continuous, so we conclude that $\tau_N(A)$ is such a linear
combination for all $A \in L^\infty_r(\T^n)$.
\bigskip

\noindent {\bf Theorem 47.} {\it Regard $L^\infty(\T^n)$ as
the subalgebra of $2\pi$-periodic functions in $L^\infty(\R^n)$.
The following equalities hold:
$$\eqalign{L^\infty_r(\T^n) &= \sigma^{-1}(L^\infty(\T^n))\cr
C_r(\T^n) &= UC_r(\R^n) \cap L^\infty_r(\T^n)\cr
{\rm Lip}_r(\T^n) &= {\rm Lip}_r(\R^n) \cap L^\infty_r(\T^n)\cr
H^\infty_r(\T) &= H^\infty_r(\R) \cap L^\infty_r(\T).\cr}$$}

\noindent {\it Proof.} Since $\sigma(T_{\vec{a}}) \in L^\infty(\T^n)$
for any $\vec{a} \in \Z^n$, linearity and weak* continuity of $\sigma$
imply that $\sigma(L^\infty_r(\T^n)) \subset L^\infty(\T^n)$. Conversely,
let $A \in L^\infty_r(\R^n)$ and suppose $\sigma(A) \in L^\infty(\T^n)$.
Although we do not yet know $A \in L^\infty_r(\T^n)$, define $\tau_N(A)$
by the formula preceding the theorem, and define $\tau_N(\sigma(A))
\in L^\infty(\R^n)$ similarly. Then $\sigma(A) \in L^\infty(\T^n)$ implies
that $\tau_N(\sigma(A)) = \sigma(\tau_N(A))$ is
a linear combination of the functions
$\sigma(T_{\vec{a}}) = e^{\iota\vec{a}\cdot\vec{p}}$
with $\vec{a} \in \Z^n$ and $-N \leq a_i \leq N$ for all $i$. Applying
$\sigma^{-1}$, we infer that $\tau_N(A)$ is in the span of the corresponding
$T_{\vec{a}}$. Now $\|\tau_N(A)\| \leq \|A\|$ since $K_N \geq 0$ and
$\int_\T K_N(t)\, dt = 2\pi$, and $\sigma(\tau_N(A)) =
\tau_N(\sigma(A)) \to \sigma(A)$
weak* is classical, so we must also have $\hat{\sigma}(\tau_N(A)) \to
\hat{\sigma}(A)$ weak operator, and hence $\tau_N(A) \to A$ weak operator
since $(\tau_N(A))$ is bounded and $L^2(\R^n)$ is dense in $L^2_r(\R^n)$.
It follows that $A$ is in the weak operator closure of elements of
${\rm span}\{T_{\vec{a}}: \vec{a} \in \Z^n\}$, so $A \in L^\infty_r(\T^n)$.
This proves the first equality.
\medskip

The forward containment of the second equality follows from the fact that
$C_r(\T^n) \subset AP_r(\R^n) \subset UC_r(\R^n)$ (Proposition 28).
For the reverse containment, let $A \in UC_r(\R^n) \cap L^\infty_r(\T^n)$.
By the remark preceding the theorem, it will suffice to show that
$\tau_N(A) \to A$ in norm. Observe that
$$A = {{1}\over{(2\pi)^n}} \int_{\T^n} A K_N^n(\vec{t}\,)\, d\vec{t}$$
since $\int_\T K_N(t)\, dt = 2\pi$. Thus
$$A - \tau_N(A) =
\int_{\T^n} (A - \gamma_{\vec{t}}(A)) K_N^n(\vec{t}\,)\, d\vec{t}.$$
This integral goes to zero in norm as $N \to \infty$ because the
map $\vec{t} \mapsto A - \gamma_{\vec{t}}(A)$ is norm continuous
and is zero when $\vec{t} = 0$, while for any $\epsilon > 0$ we have
$\int_{|\vec{t}| \geq \epsilon} K_N^n(\vec{t}\,)\, d\vec{t} \to 0$
as $N \to \infty$. Thus $\tau_N(A) \to A$ in norm.
\medskip

The third equality is trivial, as is the forward containment of
the fourth. For the reverse containment, let $A \in
H^\infty_r(\R) \cap L^\infty_r(\T)$. Then $\sigma(A)
\in H^\infty(\T)$ and $\sigma(\tau_N(A)) = \tau_N(\sigma(A))$ is a linear
combination of $\sigma(T_0), \ldots, \sigma(T_N)$. Thus $\tau_N(A)$ is
a linear combination of $T_0, \ldots, T_N$, and $\tau_N(A) \to A$
weak operator as in an earlier part of the proof, so we conclude
that $A \in H^\infty_r(\T)$.\hfill\hal
\bigskip

We can prove $C_r(\T^n) = C_{b,r}(\R^n) \cap L^\infty_r(\T^n)$ by
a somewhat more complicated argument involving passage to
$l^2(\Z^n, 1/E(\vec{k})) \otimes L^2([0,1]^n)$. The point is that
multiplication by $K_f$, for $f$ sufficiently close to a delta function
at the origin, is an isometry on $L^\infty_r(\T^n)$ with the norm
it receives in this representation. Then we use the fact that
$A \in C_{b,r}(\R^n)$ implies $AK_f \in C_{0,r}(\R^n)$, which
implies norm continuity of $\gamma_{\vec{t}}(AK_f)$; together
with $\gamma_{\vec{t}}(AK_f) \approx \gamma_{\vec{t}}(A)K_f$
for small $|\vec{t}|$, this allows us to conclude that
$\gamma_{\vec{t}}(A) \to A$ in norm as $\vec{t} \to 0$.
\bigskip

\noindent {\bf Corollary 48.} {\it Let $X_r$ be any of the spaces
in Definition 45 and let $X$ be its classical analog. Then
$\sigma(X_r) \subset X$.}
\bigskip

\noindent {\bf Corollary 49.} {\it The images under $\sigma$ of all
of the spaces in Definition 45, except $A_r(\T)$ and $H^\infty_r(\T)$,
are self-adjoint subalgebras of $L^\infty(\T^n)$.}
\bigskip

\noindent {\bf Corollary 50.} {\it We have
$$C^1_r(\T^n) \subset {\rm Lip}_r(\T^n) \subset C_r(\T^n).$$}

(For the first containment, use the comment following Proposition 29.)
\bigskip

\noindent {\bf Corollary 51.} {\it Let $S \subset \R^n$ be measurable
and periodic. Then $\I(S) \cap L^\infty_r(\T^n)$ is a weak* closed ideal
of $L^\infty_r(\T^n)$ and $\I(S) \cap C_r(\T^n)$ is a norm closed ideal
of $C_r(\T^n)$.}
\bigskip

Our last three results are not corollaries of Theorem 47. The first two
can be proven directly using the techniques of Corollary 37 and Theorem
38 in the $l^2(\Z^n, 1/E(\vec{k}))$ model of $L^\infty_r(\T^n)$. The
proof of the third resembles the proof of Theorem 40 in outline, but
is easier. The key step is to show that any complex homomorphism
$\omega: C_r(\T^n) \to \C$ takes $T_{\vec{e}_i}$ into the unit circle
($1 \leq i \leq n$). This follows from the fact that $\|T_{\vec{e}_i}^k\|
= \|T_{k\vec{e}_i}\| \sim \sqrt{k}$: then $|\omega(T_{\vec{e}_i})^k|
= |\omega(T_{\vec{e}_i}^k)| = O(\sqrt{k})$ implies $|\omega(T_{\vec{e}_i})|
\leq 1$, and the same argument applied to $T_{-\vec{e}_i} = T_{\vec{e}_i}^{-1}$
shows the reverse inequality.
\bigskip

\noindent {\bf Proposition 52.} {\it Let $d > (n+1)/2$. Then
$$C^d(\T^n) \subset \sigma(C_r(\T^n)).$$}

\noindent {\bf Proposition 53.} {\it Let $d > (n+1)/2$. Then
$\sigma(C_r^d(\T^n))$ is properly contained in
$C^d(\T^n) \cap \sigma(C_r(\T^n))$.}
\bigskip

\noindent {\bf Proposition 54.} {\it The maximal ideal space of
$C_r(\T^n)$ can be identified with $\T^n$, and the symbol map
$\sigma: C_r(\T^n) \to C(\T^n)$ can be identified with the Gelfand
transform.}
\bigskip

As a closing remark, we consider the possibility of using other weight
functions on $\R^n$ besides $1/E(\vec{p})$. For reasons explained in
Section 1, this is the only weight that interests us, but our methods
and results could apply to other weights $w(\vec{p})$ as well. What
properties of the function $1/E(\vec{p})$ were used?
\medskip

Remark 14 (a) specifically requires $w(\vec{p}) = 1/E(\vec{p})$, and
the asymptotic approximation $w(\vec{p}) \sim |\vec{p}|^{-1}$ for large
$|\vec{p}|$ is used in Remark 14 (c), Remark 23 (b), Lemma 36, Corollary 37,
Theorem 38, and Remark 46 (a). However, the bulk of the paper is valid
more broadly. The most important conditions are that $w(\vec{p})$ be
bounded (needed so that $L^2(\R^n) \subset L^2_r(\R^n)$) and that
$w(\vec{p} + \vec{a})/w(\vec{p}) \to 1$ as $\vec{p} \to \infty$, for
all $\vec{a}$ (needed for boundedness of $T_{\vec{a}}$ and $K_f$, and
in Theorem 16 and Lemma 32).
\medskip

Various other properties of $1/E(\vec{p})$ appear sporadically.
Proposition 7 requires that the weight not be integrable, Theorems
10, 11, and 12 require that it not be constant, Remark 14 (b) requires
that $w$ be a function of $|\vec{p}|$, Proposition 22 requires that
$w(-\vec{p})/w(\vec{p})$ be bounded, and Remark 23 (a) requires
that $w(\vec{p}) \to 0$ as $\vec{p} \to \infty$.
\bigskip
\bigskip

\noindent {\bf Acknowledgement}
\bigskip

The author wishes to acknowledge helpful discussions with Allan Donsig,
Nets Katz, John McCarthy, and David Pitts.
\bigskip
\bigskip

[1] W.\ Arveson, {\it An Invitation to C*-Algebras}, Springer-Verlag,
New York, 1976.
\medskip

[2] G.\ Baym, {\it Lectures on Quantum Mechanics}, W.\ A.\ Benjamin,
New York, 1969.
\medskip

[3] A.\ J.\ Bracken and G.\ F.\ Melloy, Localizing the relativistic
electron, {\it J.\ Phys.\ A: Math.\ Gen.\ \bf 32} (1999), 6127-6139.
\medskip

[4] O.\ Bratteli and D.\ W.\ Robinson, {\it Operator Algebras and
Quantum Statistical Mechanics 1: C*- and W*- Algebras, Symmetry Groups,
Decomposition of States} (second edition), Springer-Verlag, New York,
1987.
\medskip

[5] J.\ B.\ Conway, {\it A Course in Functional Analysis}, Springer-Verlag,
New York (1985).
\medskip

[6] P.\ A.\ Fillmore, {\it A User's Guide to Operator Algebras},
Wiley-Interscience, New York, 1996.
\medskip

[7] W.\ Greiner, {\it Relativistic Quantum Mechanics: Wave Equations},
Springer, New York, 2000.
\medskip

[8] R.\ Haag, {\it Local Quantum Physics: Fields, Particles, Algebras},
Springer-Verlag, New York, 1992.
\medskip

[9] G.\ C.\ Hegerfeldt, Remark on causality and particle localization,
{\it Phys.\ Rev.\ D \bf 10} (1974), 3320-3321.
\medskip

[10] E.\ M.\ Henley and W.\ Thirring, {\it Elementary Quantum Field
Theory}, McGraw-Hill, New York, 1962.
\medskip

[11] B.\ R.\ Holstein, {\it Topics in Advanced Quantum Mechanics},
Addison-Wesley, Redwood City, 1992.
\medskip

[12] J.\ A.\ Eisele, {\it Modern Quantum Mechanics, with Applications to
Elementary Particle Physics; An Introduction to Contemporary Physical
Thinking}, Wiley-Interscience, New York, 1969.
\medskip

[13] A.\ J.\ K\'alnay, The localization problem, in {\it Problems in
the Foundations of Physics} (M.\ Bunge, ed.), Springer-Verlag,
New York, 1971, pp.\ 93-110.
\medskip

[14] R.\ H.\ Landau, {\it Quantum Mechanics II: A Second Course in
Quantum Theory}, Wiley, New York, 1996.
\medskip

[15] T.\ D.\ Newton and E.\ P.\ Wigner, Localized states for elementary
systems, {\it Rev.\ Modern Phys.\ \bf 21} (1949), 400-406.
\medskip

[16] T.\ O.\ Philips, Lorentz invariant localized states, {\it Phys.\
Rev.\ \bf 136} (1964), B893-B896.
\medskip

[17] W.\ Rudin, {\it Functional Analysis} (second edition), McGraw-Hill,
New York, 1991.
\medskip

[18] P.\ Strange, {\it Relativistic Quantum Mechanics: With Applications
in Condensed Matter and Atomic Physics}, Cambridge University Press,
Cambridge, 1998.
\medskip

[19] P.\ Teller, {\it An Interpretive Introduction to Quantum Field Theory},
Princeton University Press, Princeton, NJ, 1995.
\medskip

[20] V.\ S.\ Varadarajan, {\it Geometry of Quantum Theory} (second edition),
Springer-Verlag, New York, 1985.
\medskip

[21] N.\ Weaver, {\it Mathematical Quantization}, Chapman \& Hall/CRC,
New York, 2001.
\medskip

[22] A.\ S.\ Wightman, On the localizability of quantum mechanical systems,
{\it Rev.\ Modern Phys.\ \bf 34} (1962), S45-S72.
\bigskip
\bigskip

\noindent Math Dept.

\noindent Washington University

\noindent St.\ Louis, MO 63130 USA

\noindent nweaver@math.wustl.edu
\end